\documentclass[10pt]{amsart}
\usepackage{rotating}
\usepackage{mathcomp}
\usepackage{amssymb}
\newcommand{\C}{{\Bbb C}}
\newcommand{\Z}{{\Bbb Z}}

\newcommand{\R}{{\Bbb R}}

\newcommand{\F}{{\Bbb F}}
\newcommand{\Polys}{{\mbox{\rm Polys}}}
\newcommand{\MPolys}{{\mbox{\rm MPolys}}}
\newcommand{\bbP}{{\Bbb P}}
\newcommand{\Q}{{\Bbb Q}}

\newcommand{\Gal}{\mbox{Gal}}

\newcommand{\cP}{P}

\newcommand{\disc}{\mbox{disc}}

\newcommand{\Spec}{\mbox{Spec}}
\newcommand{\Conf}{\mbox{Conf}}

\numberwithin{equation}{section}
\numberwithin{table}{section}
\numberwithin{figure}{section}
\newtheorem{Theorem}{Theorem}[section]
\newtheorem{Proposition}{Proposition}[section]
\newtheorem{ReductionBound}{Reduction Bound}[section]
\newtheorem{PullbackConstruction}{Pullback Construction}[section]

\newtheorem{Definition}{Definition}[section]

\allowdisplaybreaks
\title{Polynomials with prescribed bad primes}
\author{David P.\ Roberts}
\address{Division of Science and Mathematics, University of
  Minnesota Morris, Morris, MN 56267, USA}
\email{roberts@morris.umn.edu}

\begin{document}

\sloppy

\begin{abstract}  We tabulate polynomials in $\Z[t]$ with a given 
factorization partition, bad reduction entirely within 
a given set of primes, and satisfying auxiliary conditions 
associated to $0$, $1$, and $\infty$.   We explain how
 these sets of polynomials are of particular interest
 because of their role in the construction
of nonsolvable number fields of arbitrarily large degree and
 bounded ramification.  Finally we discuss the
 similar but technically more complicated
  tabulation problem corresponding to removing
  the auxiliary conditions.
\end{abstract}

\maketitle

\section{Introduction}  
\label{intro}

\subsection{Overview} For $\cP = \{p_1,\dots,p_r\}$ a finite set of primes, let $P^*$ be the set of integers
of the form $\pm p_1^{e_1} \cdots p_r^{e_r}$.  We say that a polynomial in $\Z[t]$ is normalized
if its leading coefficient $s(\infty)$ is positive and the greatest common divisor of its coefficients is $1$.  \begin{Definition}
\label{maindef}
For $\kappa$ a partition, 
 $\Polys_\kappa(\cP)$ is the set
of normalized polynomials $s(t) \in \Z[t]$ satisfying 
\begin{description}
\item[1] The degrees of the irreducible factors of $s(t)$ form the partition $\kappa$;
\item[2] The discriminant $\mbox{{\rm Disc}}(s)$ and the values $s(0)$, $s(1)$, $s(\infty)$
are all in $P^*$.   
\end{description}
\end{Definition}
\noindent   The results of this paper identify 
many $\Polys_\kappa(\cP)$ completely and show that others 
are large.  

A sample theoretical result and some computational results within it 
give a first sense of the content of this paper.  
The theoretical result is an algorithm to 
determine $\Polys_{3^c 2^b 1^a}(P)$ given the
set of all $j$-invariants of elliptic curves with bad reduction
within $P \cup \{2,3\}$.     The computational result uses Coghlan's
determination \cite{Co} of the eighty-three $j$-invariants for $P = \{2,3\}$
as input.  Carrying out the algorithm gives  
$\Polys_{3^c 2^b 1^a}({\{2,3\}})$ for all $(a,b,c) \in \Z_{\geq 0}^3$.    
The largest cardinality arising is $|\Polys_{3^4 1}(\{2,3\})| = 180,822$.  
The largest degree $3c+2b+a$ coming from a nonempty set of polynomials is $35$,
arising uniquely from $|\Polys_{3^{11}2}(\{2,3\})|=2$.  One of the 
two elements of $\Polys_{3^{11}2}(\{2,3\})$ is 
\begin{eqnarray}
\nonumber s(t) & = &
\nonumber   
  \left(t^3-2\right) \left(t^3+3 t^2-3 t+1\right) \left(2 t^3-6 t^2+6 t-1\right) \cdot \\
\label{big23}    &&  \left(t^3-3 t+4\right) \left(2 t^3+3 t-1\right) \left(4   t^3-9 t^2+6 t-2\right) \cdot \\
\nonumber    && \left(t^3-3 t^2+6 t-2\right) \left(2 t^3-3 t+2\right) \left(2 t^3-3 t^2-1\right) \cdot \\
\nonumber    && \left(t^3-3 t+1\right) \cdot  \left( t^3-3 t^2+1\right)  \cdot \left(t^2-t+1\right). 
 \end{eqnarray}
 The other one is $t^{35} s(1/t)$, and both polynomials have discriminant $2^{105} 3^{533}$.  
 
   Our primary motivation
is external, as polynomials in $\Polys_\kappa(P)$
 are
used in the construction of two types of 
nonsolvable number fields of arbitrarily large degree
 and bounded ramification.   
 Katz number fields \cite{RABC}, \cite{RM05}
 have Lie-type Galois groups and 
 the least ramified examples tend 
 to have two ramifying primes.  
 Hurwitz number fields 
 \cite{RHNF, RV} 
typically have alternating or symmetric 
 Galois groups and the least ramified
 examples tend to have 
 three ramifying primes.

    The natural problem corresponding to our title 
involves suitably tabulating polynomials when the conditions 
$s(0)$, $s(1)$, $s(\infty) \in P^*$ are removed. The special case we
pursue here is more elementary but has much of the
character of the general problem. 
   The full problem is 
briefly discussed at the end of this paper.  

\subsection{Three steps and three regimes}
\label{33}
Constructing all elements of $\Polys_\kappa(\cP)$ in general is naturally a three-step process. 
Step 1 is to identify the
set $NF_d(P)$ of isomorphism classes of degree $d$ number fields 
ramified within $P$, for each $d$ appearing in $\kappa$.  
For many $(d,P)$ this complete list is available at \cite{Database}.
 Step 2 is to get the contribution $\Polys^K_d(\cP)$ of 
 each $K \in NF_d(P)$ to $\Polys_d(\cP)$ by inspecting
 the finite set of exceptional $P$-units in $K$. 
 We expect an algorithm finding these units 
to appear in standard 
 software shortly, generalizing the implementation
 in {\em Magma} \cite{Mag} for the case $\cP = \emptyset$.   
Step~3 is to extract those products of the irreducible
polynomials which are in $\Polys_{\kappa}(P)$.  
This last step is essentially bookkeeping, but
nonetheless presents difficulties as $\Polys_\kappa(P)$ can
be very large even when all the relevant
$\Polys_d(P)$ are relatively small.  

One can informally distinguish three regimes as follows.  
For suitably small $(\kappa,\cP)$, one can ask for the provably complete
list of all elements in $\Polys_\kappa(\cP)$.   For intermediate
$(\kappa,\cP)$, one can seek lists which seem likely to be
complete.  For large $(\kappa,\cP)$, one can
seek systematic methods of constructing interesting
 elements of $\Polys_\kappa(\cP)$.  
 We present  results here in all 
 three regimes.

\subsection{Content of the sections}
Section~\ref{graph} consist of preliminaries, with a focus 
on carrying out Step~3 by 
interpreting polynomials in $\Polys_\kappa(\cP)$ as 
cliques in a graph $\Gamma(\cP)$.  
Sections~\ref{1polys}, \ref{2polys}, and \ref{3polys} 
are in the first regime and are similar to each other 
in structure.     They present
general results corresponding to partitions $\kappa$ 
of the form $1^a$, $2^b1^a$, and $3^c2^b1^a$ 
respectively.  In these results, Steps~1 and 2 are carried out together
by techniques particular to $d \leq 3$ involving ABC triples.
 As illustrations of the generalities, these sections
completely identify all $\Polys_{1^a}(\{2,3,5,7\})$,
 $\Polys_{2^b 1^a}(\{2,3,5\})$,
 and the above-discussed $\Polys_{3^c 2^b 1^a}(\{2,3\})$.

  Section~\ref{genpolys} is in the second regime
  and follows the three-step approach.  
  To illustrate the general method, this
  section takes $P = \{2\}$ so
  that $NF_d(\{2\})$ is known
  to be empty for $d \in \{3,5,6,7\}$.  
 It identifies
 all $\Polys_{4^d 2^b 1^a}(\{2\})$, assuming
 the identification of $\Polys_{4}(\{2\})$
 is correct.  
Because of the increase in allowed $\kappa$ in
Sections~\ref{1polys}-\ref{genpolys}, our considerations
become conceptually more complicated.  
Because of the simultaneous decrease in $\cP$, 
our computational examples remain at approximately the
same level of complexity.
 Section~\ref{largedegree} is in the third regime.  
 It shows that some $\Polys_\kappa(\cP)$ are large because of products of
cyclotomic polynomials while  others are large because 
of polynomials related to fractals. 

  Section~\ref{specialization} 
sketches the applications to number field construction.   
 Our  presentation gives
a feel for how the $\Polys_\kappa(P)$ enter by 
presenting one family of examples from the Katz setting 
and one family from the Hurwitz setting.  
Section~\ref{conclude} concludes the paper
 by discussing promising directions for future work,
 with a focus on moving into the more
 general setting where the auxiliary
 conditions on $s(0)$, $s(1)$, and $s(\infty)$ are 
 removed.

 \subsection{Acknowledgements}
 We thank Frits Beukers, Michael Bennett, John Cremona, John Jones, and Akshay Venkatesh for
 conversations helpful to this paper.  We thank the Simons
 Foundation for research support through 
  grant   \#209472.

\section{Preliminaries}
\label{graph}

\subsection{Sets related to $\Polys_\kappa(\cP)$.}  It is convenient to consider disjoint unions of
$\Polys_\kappa(\cP)$ over varying $\kappa$ as follows:
\[
\begin{array}{rlcrl} 
\multicolumn{2}{c}{\mbox{Full sets}}  & &\multicolumn{2}{c}{\mbox{Finite subsets}}   \\
\cline{1-2} \cline{4-5} 
\Polys(\cP)  &\!\! =  {\displaystyle \coprod_{\kappa} \Polys_\kappa(\cP),} & &
\Polys(\cP)^f & \!\!=  {\displaystyle \coprod_{{\rm max}(\kappa) \leq f}\Polys_\kappa(\cP),} \\
\Polys(\cP)_\ell & \!\!=  {\displaystyle \coprod_{{\rm length}(\kappa) = \ell}\Polys_\kappa(\cP),} & \;\;\;\;\; & 
\Polys(\cP)_\ell^f & \!\! =  \Polys(\cP)^f \cap \Polys(\cP)_\ell.
\end{array}
\]
Thus $\Polys(\cP)$ is the set of all polynomials under 
study for a given $\cP$.   It and the subsets
$\Polys(\cP)_\ell$ are always infinite for 
any $\cP \neq \emptyset$ and $\ell \geq 1$, 
as discussed further in Section~\ref{largedegree}.

We say that a polynomial is {\em $f$-split} if all its irreducible factors
have degree at most $f$.  From more general theorems cited 
in Section~\ref{conclude},  the sets $\Polys(\cP)^f$ and 
thus $\Polys(\cP)^f_\ell$ of $f$-split polynomials 
are always finite.  To focus just on degree and suppress reference
to the factorization partition, another convenient finite set is
$\Polys_{[k]}(P) = \coprod_{\kappa \vdash k} \Polys_\kappa(P).$  

\subsection{Compatibility} The study of $\Polys(\cP)$ reduces to a great extent
to the study of $\Polys(\cP)_1$ as follows.  Let $s_1$, \dots, $s_\ell$ 
be in $\Polys(\cP)_1$, thus irreducible normalized polynomials in 
$\Z[t]$, with discriminants $D_i$ and 
values $s_i(0)$, $s_i(1)$, $s_i(\infty)$ all in $\cP^*$.   
The product $s(t) = s_1(t) \cdots s_\ell(t)$ 
certainly satisfies $s(0)$, $s(1)$, $s(\infty) \in \cP^*$.  
Its discriminant is given by the product formula
\[
D = \left( \prod_{i=1}^\ell D_i \right) \left( \prod_{i<j} R^2_{ij} \right),
\]
where $R_{ij}$ is the resultant $\mbox{Res}(s_i,s_j) \in \Z$.  
In general, we say that two polynomials $u$ and $v$ in $\Polys(\cP)$ are 
{\em compatible} if $\mbox{Res}(u,v) \in \cP^*$.  
Thus $s \in \Polys(\cP)_\ell$ if and only if
its $\ell$ irreducible factors are pairwise compatible.  

\subsection{Graph-theoretic interpretation}  To exploit the notion of compatibility,
we think in terms of a graph $\Gamma(\cP)$ as follows.   
The vertex set of $\Gamma(\cP)$ is 
$\Polys(\cP)_1$.    
If a vertex corresponds to a degree $d$ polynomial, 
we say it has degree $d$.
The edge-set
of $\Gamma(\cP)$ is 
$\Polys(\cP)_2$,
 with an edge $s_1s_2$ having
endpoints $s_1$ and $s_2$.  Thus edges 
are placed between compatible irreducible polynomials.
In general, a polynomial in  $\Polys(\cP)_\ell$ 
is identified with a clique in $\Gamma(\cP)$ of size $\ell$, 
meaning a complete subgraph on $\ell$ vertices.   
For similar use of graph-theoretic language in
contexts like ours, see e.g.\  \cite{LN}.  

When restricting attention to $f$-split polynomials, we likewise
think in terms of the corresponding graph $\Gamma(\cP)^f$.  
This graph is now finite, with vertex set $ \Polys(P)^f_1$, 
edge set $ \Polys(P)^f_2$, 
and cliques of size $\ell$ 
corresponding to elements of $\Polys(\cP)_\ell^f$.  
Figure~\ref{twograph}, discussed in more detail in \S\ref{2graph} below, draws
$\Gamma(\{2\})^2$.

\subsection{Packing points into the projective line}  Our problem of identifying 
$\Polys_\kappa(\cP)$ 
can be understood in geometric language as follows.  For each
prime $p$, let $\overline{\F}_p$ be an algebraic closure
of $\F_p$.  For any prime power $p^f$, let 
$\F_{p^f}$ be the subfield of $\overline{\F}_p$ having 
$p^f$ elements.   For any field $F$, let $\bbP^1(F) = F \cup \{\infty\}$
be the corresponding projective line.  

Let $s(t) \in \Polys(\cP)$ have degree $k$.  Denote its set of complex roots by $Z$,
so that $|Z|=k$.    Let $\widehat{Z} = Z \cup \{0,1,\infty\} \subset \bbP^1(\C)$.   
For any prime $p$, similarly let $Z_p$ be the root-set of $s(t)$ in $\overline{\F}_p$ 
and $\widehat{Z}_p = Z \cup \{0,1,\infty\} \subset \bbP^1(\overline{\F}_p)$. 

Let $\overline{\Q} \subset \C$ be the field of algebraic numbers.     
Via roots, our $\Polys(\cP)$ is in bijection with the set of finite subsets  
$\widehat{Z} \subset \bbP^1(\overline{\Q})$ which are $\Gal(\overline{\Q}/\Q)$-stable,
contain $\{0,1,\infty\}$, and have good reduction outside of $\cP$ in the sense that
the reduced sets $\widehat{Z}_p$ have the same
size as $\widehat{Z}$ for $p$ a prime not in $\cP$.      

If $s(t) \in \Polys(\cP)^f$ then the set $\widehat{Z}_p$ lies in the finite
set
$
\cup_{d \leq f} \bbP^1(\F_{p^d}).
$
The order of this set for $f = 1$, $2$, $3$, and $4$ is respectively 
$p+1$, $p^2+1$, $p^3 + p^2 - p+1$, $p^4 + p^3 - p + 1$.   
One has the following trivial bound, which we highlight because of its importance:
\begin{ReductionBound} 
\label{trivbound} 
A polynomial $s(t) \in \Polys(\cP)^f$ has degree 
at most 
\[
N(p,f) = |\bigcup_{d \leq f} \bbP^1(\F_{p^d})| - 3,
\]
where $p$ is the 
smallest prime not in $\cP$.  
\end{ReductionBound}
\noindent The room available for packing points increases polynomially
with the first good prime  $p$ and exponentially with the degree cutoff $f$:
\[
 \begin{array}{r|rrrrl}
             N(p,f)  &    1 & 2 & 3 & 4 &   \\
                   \cline{1-5}
             2 &      0 & 2 & 8 & 20 \\
             3 &      {\mathit{1}} & {\mathit{7}} & 31 & {\bf 103} \\
            5  &      3 & 23 & {\bf 143} & 743 \\
            7 &       5 & {\bf 47} & 383 & 2735 \\
            11 &    {\bf 9} & 119 & 1439 & 15959 & \!\!\!\!\!\! . \!\!\! \\
                  \end{array}
\]
The italicized entries are relevant to Figure~\ref{twograph} 
where both upper bounds are achieved.  The boldface entries ascending
 to the right correspond to Sections~\ref{1polys}, \ref{2polys}, \ref{3polys}, \ref{genpolys}
respectively, with the bound being obtained only in the first case.

\subsection{$S_3$-symmetry}  
\label{S3} If $s(t) \in \Polys_{\kappa}(\cP)$ has degree $k$, then its properly signed
transforms
\begin{eqnarray*}
s_{(01)}(t) & = & \pm s(1-t), \\
s_{(0\infty)}(t) & = & \pm t^k s(1/t), 
\end{eqnarray*}
are also elements in $\Polys_{\kappa}(\cP)$.    The two displayed
transformations generate a six-element group $S_3$ which 
acts on each $\Polys_{\kappa}(\cP)$.  
Our notation captures that these transformations arise from permuting
the special points $0$, $1$, and $\infty$ arbitrarily. 

\subsection{The graph $\Gamma(\{2\})^2$.}
\label{2graph}
\begin{figure}[htb]
\begin{tabular}{cc}
\includegraphics[width=3in]{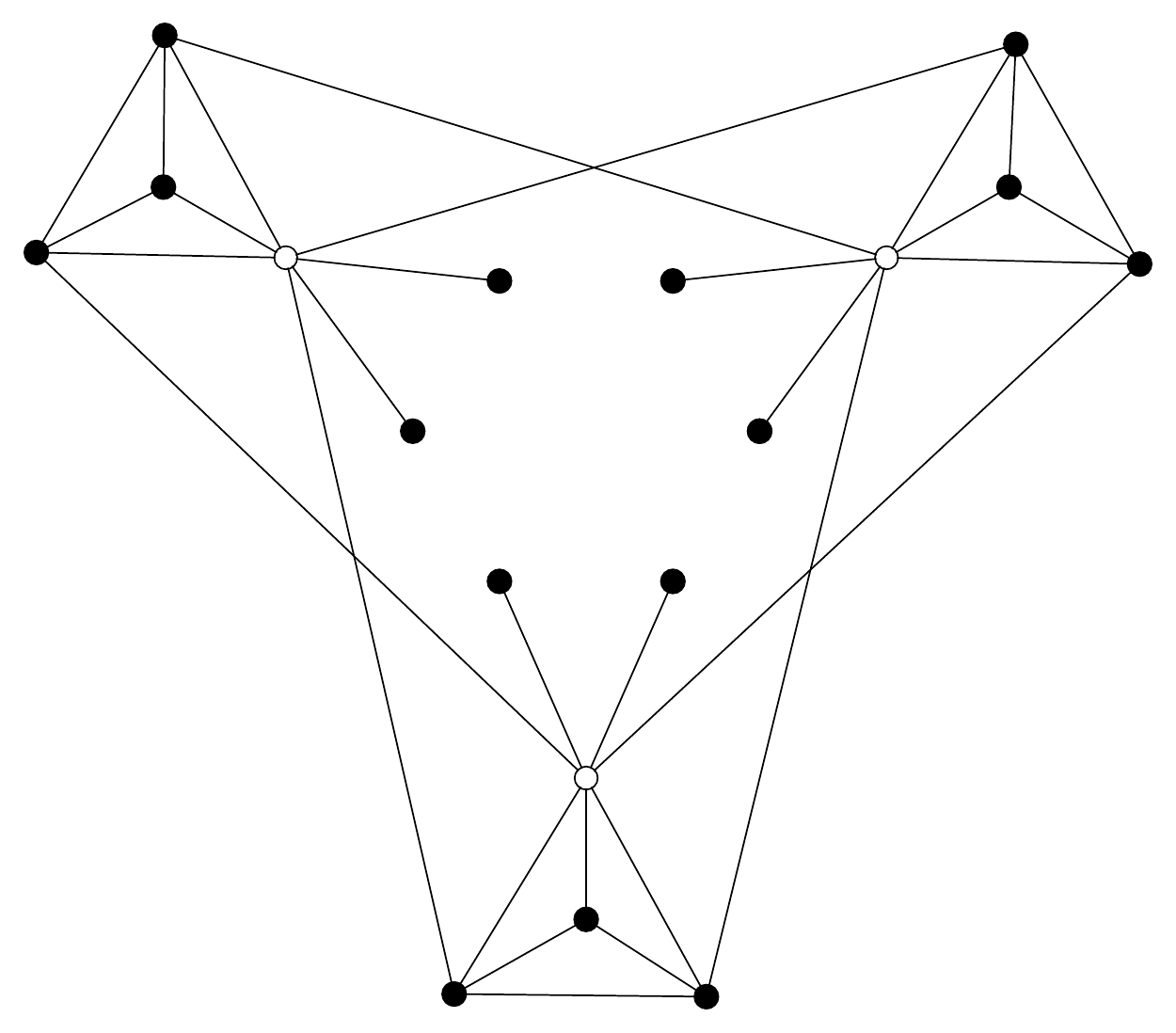} & $ \!\!\!\!\!\!\!\!\!\! \begin{array}{ccc}
                   t^2-8 t+8, \!\!\!\!\!\!  \!\!\!\!\!\! \!\!\!\!\!\!   \!\!\!\!\!\! &  & \!\!\!\!\!\!  \!\!\!\!\!\!  \!\!\!\!\!\! \!\!\!\!\!\!  t^2+4 t-4 \\
                   &&\\
                    & t-2 &  \\
                    &&\\
                    & t^2-2 t+2 &  \\
                    &&\\
                   t^2-2,  \!\!\!\!\!\!  &  &  \!\!\!\!\!\!  t^2-4 t+2 \\
                   &\\
                   &\\
                   &\\
                   &\\
                   &\\
                   &\\
                   &\\
                  \end{array}
$
\end{tabular}
\vspace{-1in}
\caption{\label{twograph}  The graph $\Gamma(\{2\})^2$.
The polynomials represented by the vertices in the lower third of the graph 
are indicated.  
}
\end{figure}

    Figure~\ref{twograph}
 gives a simple example illustrating many of our considerations so far.  
 The three white vertices are the polynomials
 in $\Polys_1(\{2\})$ and the subgraph $\Gamma(\{2\})^1$ consists
 of three isolated points.  The fifteen black vertices are the elements of
$\Polys_2(\{2\})$.  That the drawn graph is indeed all
of $\Gamma(\{2\})^2$ is a special case of the completeness
results cited in Section~\ref{2polys}.

     The sets $\Polys_{2^b 1^a}(\{2\})$ can all be read off of 
Figure~\ref{twograph}, and have sizes
\begin{equation}
\label{littletab2}
\begin{array}{r|rrl}
b & a=0  & a=1   \\
\cline{1-3}
0 & 1 & 3  \\
1 & 15 & 21 \\
2 & 9 & 9 \\
3 & 3 & 3 & \!\!\!\!\!\! . \!\!\!
\end{array}
\end{equation}
For example, the clique formed by the four lowest vertices gives the element
\[
s(t) = (t^2-2t+2)(t^2-2)(t^2-4t+2)(t-2)
\]
 of $\Polys_{2^31}(\{2\})$.  The polynomial $s(t)$ and its transforms by $S_3$ give the bottom right
 $3$ of \eqref{littletab2}.    The graph-theoretic deductions in Sections~\ref{1polys}-\ref{genpolys} 
 are conceptually no different from the visual inspection of Figure~\ref{twograph} needed to produce 
 \eqref{littletab2}.  However the graphs involved are much larger and the
 passage from graphs to cliques is incorporated into our
 programs as described next.

\subsection{Step 3 of the process}
\label{twostep}
    To compute a graph $\Gamma(\cP)^f$ and all associated sets $\Polys_\kappa(\cP)$,
 the first two steps as described in \S\ref{33} yield the 
 vertex   set $\Polys(\cP)^f_1$.  
 Step 3, passing from the vertex set
  to the entire graph,
  is then done as follows.  For each vertex $s_1(t) \in \Polys(\cP)^f_1$  we compute resultants
 and determine its set $N_{s_1(t)}$ of lesser neighbors with respect
 to some ordering.  The edge set $\Polys(\cP)^f_2$ is then all
 $s_1(t)s_2(t)$ with $s_2(t) \in N_{s_1(t)}$.   One continues 
 inductively, with $\Polys(\cP)^f_\ell$ being the 
 set of $s_1(t) \cdots s_\ell(t)$ with $s_\ell(t) \in \cap_{i=1}^{\ell-1} N_{s_i(t)}$.

\subsection{Monic variants}
\label{monic}
   If $s(t) \in \Z[t]$ is a normalized polynomial then $s(t)/s(\infty) \in \Q[x]$
is a monic polynomial.  It is often technically more convenient to work
with monic rather than normalized polynomials.   
Accordingly, we let $\MPolys(P)$ be the set of monic polynomials
$s(t)/s(\infty)$ with $s(t) \in \Polys(P)$.    So elements
of $\MPolys(P)$ lie in $\Z^P[t]$, where
 $\Z^P$ is the ring of rational numbers with 
 denominators in $P^*$.    As a general
 rule, we keep the focus on $\Polys(P)$,
 switching temporarily to the very mild variant $\MPolys(P)$ 
 only when it is truly preferable.

\section{1-split  polynomials}
\label{1polys}
  This section describes how one determines the sets $\Polys_{1^{a}}(\cP)$,
  We illustrate the  procedure by determining $\Polys_{1^{a}}(\{2,3,5,7\})$
  for all $a \in \Z_{\geq 0}$. 
 
 \subsection{Vertices via ABC triples}  Step~1 from the introduction is trivial, since the only degree one number field is $\Q$.  
  Step~2 is to determine the 
 polynomials 
 which lie in the vertex set $\Polys_1(P)$ of the 
 graph $\Gamma(P)^1$.    To make the $S_3$-symmetry of \S\ref{S3}
 completely evident it is convenient to work with ABC triples.
\begin{Definition}
\label{abc}  For a rational number $u \neq 0,1$, let $A$, $B$, and $C$
be the unique pairwise relatively prime integers with $u = -A/C$,
 $A+B+C=0$, and $ABC<0$.    For a set of primes $P$, 
 the set $T_{\infty,\infty,\infty}(\Z^P)$ is the set of $u$ such that 
 $A$, $B$, and $C$ are in $P^*$.  
 \end{Definition}
\noindent The notation $T_{\infty,\infty,\infty}(\Z^P)$ is a specialization of the
 general notation $T_{p,q,r}(\Z^P)$ of \cite{RABC}, and we will use
 other special cases in the next two sections.  
 The action of $S_3$ on ABC triples by permutations
 corresponds to an action of $S_3$ on the projective $u$-line
 by fractional linear transformations, with $(AB)$ corresponding
 to $u \mapsto 1-u$ and $(AC)$ to $u \mapsto 1/u$.   Using
 the alternative monic language of \S\ref{monic}, one
 has 
 \[
 \MPolys_1(P) = \{t - u\}_{u \in T_{\infty,\infty,\infty}(\Z^P)}.
 \]
 This very simple parametrization is a prototype 
 for the more complicated parametrizations given
 in Theorems~\ref{thm2111} and \ref{thm3111}.
 
The set $T_{\infty,\infty,\infty}(\Z^P)$ is empty if $2 \not \in P$ by Reduction Bound~\ref{trivbound}.    
Otherwise $\{-1,1/2,2\}$ is a
three-element $S_3$-orbit and all other $S_3$-orbits have 
size six.   Elements of $T_{\infty,\infty,\infty}(\Z^P)$ can be found by computer
searches:  to get all those with 
$\mbox{height}(u) := \max(|A|,|C|)$ less than
a certain cutoff, one searches over candidate $(A,C)$
and selects those for which $B=-A-C$ is
also in $P^*$. 

      In the case $P = \{2,3,5,7\}$, a search up to 
 height $10^9$ took ten seconds and yielded 
 $375 = 3 + 6 \cdot 62$ elements.    The eighteen of largest
height come from the ABC triples 
\begin{eqnarray*}
(1,4374,-4375) & = & (1,2^13^7,-5^47), \\
(1,2400,-2401) & = & (1,2^53^15^2,-7^4),\\
(5,1024,-1029) & = & (5,2^{10},-3^1 7^3). 
\end{eqnarray*}
All the other elements have height at most 625.     The 
completeness of this list is a special case of a result of de Weger \cite[Theorem~5.4]{dW}.
This result also gives
 $|\Polys_1(\{2,3,5,7,11\})| = 1137$ 
and  $|\Polys_1(\{2,3,5,7,11,13\})| = 3267$, with largest heights
$18255$ and $1771561$ respectively.

\subsection{The sets $\Polys_{1^a}({\{2,3,5,7\}})$ }  
Tabulating cliques as described in \S\ref{twostep} has a run-time of about two minutes 
and gives the following result.
\begin{Proposition}  The nonempty sets $\Polys_{1^a}(\{2,3,5,7\})$ have size as follows:
{\renewcommand{\arraycolsep}{4pt}
\[
{
\begin{array}{r|rrrrrrrrrr}
a & 0 & 1 & 2 & 3 & 4 & 5 & 6 & 7 & 8 & 9 \\
\hline
\mbox{{\rm Size}} &1 & 375 & 9900 & 73000 & 232260 & 383712 & 356916 & 190620 & 55935 & 7425 
\end{array}
}
\]
}
\end{Proposition}
\noindent The sets involved in the next case $\{2,3,5,7,11\}$ are already much
larger, both because of the larger vertex set and from the relaxation of the
compatibility condition.

\subsection{Extremal polynomials} 
\label{extremal1} One of the $7425$ elements of $\Polys_{1^9}(\{2,3,5,7\})$ is
$s(t) = \prod_{u=2}^{10} (t-u)$.    Similarly, suppose $P$ consists of all primes strictly
less than a fixed prime $p$.  Then $s(t) = \prod_{u=2}^{p-1} (t-u)$ realizes Reduction Bound~\ref{trivbound}.

The $7425$ polynomials in $\Polys_{1^9}(\{2,3,5,7\})$ are structured into packets as follows.
Let $\prod_{i=1}^9 (t-u_i)$ be a polynomial in $\Polys_{1^9}(\{2,3,5,7\})$ 
and consider the twelve element set $\{u_1,\dots,u_9,0,1,\infty\}$.  For any
triple of distinct elements there is a unique fractional linear transformation
in $PGL_2(\Q)$ which takes these elements in order to $0$, $1$, and $\infty$.  A given 
element of  $\Polys_{1^9}(\{2,3,5,7\})$ determines $12 \cdot 11 \cdot 10/|A|$ 
elements of $\Polys_{1^9}(\{2,3,5,7\})$ in this way, with $A$ its stabilizer
subgroup.    There are in fact thirteen such packets, eight with stabilizer subgroup
$C_2$ and one each with stabilizer $C_1$, $V$, $S_3$, $D_4$ and $D_6$.
The product  $\prod_{u=2}^{10} (t-u)$ is in one of the eight packets with stabilizer $C_2$,
its nontrivial automorphism being $t \mapsto 11-t$.  As another example, the element
\[
s(t) = (t+14)(t+8)(t+5)(t+4)(t+2)(t-2)(t-4)(t-10)(t-16)
\]
represents the packet with trivial stabilizer $C_1$.  The numbers presented
are consistent via the mass-check
\[
\frac{7425}{12 \cdot 11 \cdot 10} = 
5.875 = 1 + 8 \cdot  \frac{1}{2} +  \frac{1}{4} + \frac{1}{6} + \frac{1}{8} + \frac{1}{12}.
\]
The two minute run-time cited above corresponds to a simple program which does
not exploit this type of symmetry.

\section{2-split polynomials}
\label{2polys}
    This section describes how one determines sets $\Polys_{2^b1^{a}}(\cP)$.
Without loss of generality we restrict to $\cP$ containing $2$ throughout this section.   
Assuming $\Polys_1(\cP)$ as known from the previous
section, to complete Steps 1 and 2 one needs to 
determine $\Polys_2(\cP)$ and Theorem~\ref{thm2111} gives our method.  We illustrate
the full procedure by determining all $\Polys_{2^b 1^{a}}(\{2,3,5\})$.
 
\subsection{Vertices via $ABC$ triples}  Let $T_{\infty,2,\infty}(\Z^P)$ be the 
set of rational numbers $w = -A/C$ exactly as in Definition~\ref{abc} except that
$B = -A-C$ is only required to have the form $by^2$ with $b \in  P^*$.   
For an element $w \in T_{\infty,2,\infty}(\Z^P)$, its discriminant 
class by definition is $\delta = w(1-w) \in \Q^{\times}/\Q^{\times 2}$.  
This invariant gives a decomposition
\[
T_{\infty,2,\infty}(\Z^P) = \coprod_\delta T_{\infty,2,\infty}(\Z^P)^\delta.
\]
This decomposition is used in Theorem~\ref{thm2111} below as one of two aspects of
compatibility.  

To find all $w$ in $T_{\infty,2,\infty}(\Z^P)$ up to a height bound of $H$, one searches over the
exact same set of $(A,C)$ as in the search for elements of $T_{\infty,\infty,\infty}(\Z^P)$. However now one keeps 
those $(A,C)$ where the square-free part of $B=-A-C$ is in $P^*$.   For our example,
we need the set $T_{\infty,2,\infty}(\Z^{\{2,3,5\}})$.  A one-second search
up to cutoff $H = 10^9$ found $183$ elements. 
The list consists of $-1$ and then $92$ reciprocal pairs.  
The three pairs of
largest height come from the triples
\begin{eqnarray*}
(1,-25921,25920) & = & (1,161^2, 2^6 3^4 5), \\
(9,-64009,64000) & = & (3^2, -253^2, 2^9 5^3), \\
(15625,-17161,1536) & = & (5^6, - 131^2,2^9 3).
\end{eqnarray*}
All the other elements have height at most 6561.  The completeness
of this list follows from \cite{Cr}, where the larger
set $T_{3,2,\infty}(\Z^{\{2,3,5\}})$ is calculated to have 440 elements.   The distribution according to discriminant class $\delta$ 
is quite uneven, and given
after the proof of Theorem~\ref{thm2111} 
below. 

\subsection{From the set $T_{\infty,2,\infty}(\Z^P)$ of $ABC$ triples to the set $\Polys_2(P)$ of degree two vertices in
$\Gamma(P)$}

A general quadratic polynomial in $\Q[x]$ can be written uniquely in the form
 \begin{equation}
 \label{quadformu}
 s(u_0,u_1,u_\infty; t) = u_\infty t^2 + (u_1-u_0-u_\infty) t + u_0
 \end{equation}
with $(u_0,u_1,u_\infty) \in \Q^3$.  Its discriminant is
 \begin{eqnarray}
 \Delta(u_0,u_1,u_\infty) & =  & (u_1-u_0-u_\infty)^2 - 4 u_0 u_\infty \\ 
\nonumber & = & \left(u_0^2 + u_1^2 + u_\infty^2 \right)
   -2 \left( u_0 u_1 + u_0 u_\infty + u_1 u_\infty \right).
 \end{eqnarray}
To complete an identification of the new part $\Polys_2(P)$ of the vertex set, we use
the following result, which naturally gives $\Polys_{[2]}(P) = \Polys_2(P) \coprod \Polys_{1^2}(P)$.  
 
 \begin{Theorem} 
 \label{thm2111}  Let $P$ be a finite set of primes containing $2$.   Let 
$(\delta; w_0,w_1,w_\infty)$ run over triples where $\delta \in P^*$ 
is a square-free integer and the  
 $w_i$ are in $T_{\infty,2,\infty}(\Z^P)^\delta \cup \{1\}$ satisfying
 \begin{equation}
 \label{wrel} 
\Delta(w_0,w_1,w_\infty) = - 4 w_0 w_1 w_\infty. 
\end{equation}
Then the polynomials 
\begin{equation}
\label{quadformw}
S(w_0,w_1,w_\infty) = \frac{1}{w_\infty} s(w_0,w_1,w_\infty;t)
 \end{equation}
have discriminant class $\delta$ and run over $\MPolys_{[2]}(P)$ 
\end{Theorem}

\proof   Just using that $w_0$, $w_1$, $w_\infty$ are all in $\Z^{P \times}$ one
immediately gets that $S(0)$, $S(1)$, $S(\infty)$ are all in $\Z^{P \times}$.  
Assuming further that $(w_0,w_1,w_\infty)$ satisfies \eqref{wrel}, then
the discriminant $\Delta(w_0,w_1,w_\infty)$
is also in $\Z^{P \times}$.  Thus quadratic polynomials as in 
the theorem are indeed in $\MPolys_{[2]}(P)$.  The issue
which remains is that  these polynomials form all of
$\MPolys_{[2]}(P)$.  To prove this converse direction
we start with the hypothesis that $s(u_0,u_1,u_\infty;t)/u_\infty \in \MPolys_{[2]}(P)$
and deduce that $(u_0,u_1,u_\infty)$ is proportional to a triple
$(w_0,w_1,w_\infty)$ as in the theorem.  

In general, suppose given an ordered triple of disjoint divisors $(D_2, D_{1a}, D_{1b})$ 
on the projective line $\bbP^1$ over $\Q$, of degrees $2$, $1$, and $1$ respectively.
After applying a fractional linear transformation, one can partially normalize
so that  $D_{1a} = \{0\}$, $D_{1b} = \{\infty\}$, and $D_2$ consists
of the roots of $t^2 + b t + c$ with $b$, $c \in \Q$.  
 To continue with the normalization, 
suppose $b \neq 0$.  Then one can uniquely scale so that 
still $D_{1a} = \{0\}$ and $D_{1b} = \{\infty\}$ but now 
$D_2$ consists of the roots of $t^2-t+v$ for $v = c/b^2$ in $\Q$.  
Writing $v =1/4(1-w)$, one gets that $PGL_2(\Q)$-orbits
of the initial tuple $(D_2,D_{1a},D_{1b})$ yielding $b \neq 0$ 
are in bijection with $w \in \Q - \{0,1\}$.   Moreover, the 
discriminant class in $\Q^\times/\Q^{\times 2}$ of the divisor $D_2$ is
$w(1-w)$.  Moreover, the orbit has a representative with good reduction
outside of $P$ if and only if $w \in T_{\infty,2,\infty}(\Z^P)$. 
There are infinitely
many different orbits yielding $b=0$ and we associate
all of them to $w=1$.

Let $Z$ be the roots of \eqref{quadformu}.  Then the invariants associated to 
$(Z,\{1\},\{\infty\})$, $(Z,\{\infty\},0)$, and $(Z,\{0\},\{1\})$ work out respectively to  
\[
(w_0,w_1,w_\infty)  = \frac{-\Delta(u_0,u_1,u_\infty)}{4 u_0 u_1 u_\infty} (u_0,u_1,u_\infty). 
\]
Thus any
element of $\MPolys_{[2]}(P)$ is indeed of the special form \eqref{quadformw}.  \qed

Let $\Polys^\delta_{[2]}(P)$ be the subset of 
$\Polys_{[2]}(P)$ consisting of polynomials of discriminant class $\delta$.  
Applying Theorem~\ref{thm2111} for $P = \{2,3,5\}$ to the known set
$T_{\infty,2,\infty}(\Z^P)$ gives sizes as follows:
\begin{equation*}
{\renewcommand{\arraycolsep}{2pt}
\begin{array}{r|rrrrrrrrrrrrrrrr}
\delta &    -30 & -15 & -10 & -6 & -5 & -3 & -2 & -1 & 1 & 2 & 3 & 5 & 6 & 10 & 15 &
      30 \\
      \hline
|T_{\infty,2,\infty}(\Z^P)^{\delta}| &     3 & 6 & 24 & 25 & 11 & 8 & 6 & 49 & 12 & 9 & 2 & 9 & 6 & 0 & 13 & 0 \\
|\Polys^\delta_{[2]}(\Z^P)| &  12 & 48 & 456 & 504 & 138 & 84 & 48 & 1020 & 171 & 108 & 10 & 96 & 48 & 0 &
      204 & 0 
   \end{array}
   }
\end{equation*}
Define the height of a normalized polynomial \eqref{quadformu} to be $\max(|u_0|,|u_1|,|u_\infty|)$.  
 With this definition, the height of a polynomial \eqref{quadformu} depends 
only on its $S_3$ orbit.  The three $S_3$-orbits with largest height all have height $3125=5^5$.  They are represented
by the following elements:
\begin{align*}
(w_0,w_1,w_\infty) & = \frac{(-3^7,2^7,5^5)}{5^3},    & s(t) & = 3125 t^2-810 t-2187  \\
(w_0,w_1,w_\infty) & = \frac{(-3^3,2^{11},5^5)}{2^5 3 \! \cdot \! 5}, &  s(t) & = (25 t-9) (125 t+3), \\
(w_0,w_1,w_\infty) & = \frac{(-3,2^{10} 3,5^5)}{2^6 3^2 5}, &  s(t) & = (25 t-1) (125 t+3). 
\end{align*}

\subsection{The sets $\Polys_{2^b 1^{a}}(\{2,3,5\})$}  Inductively tabulating cliques in
$\Gamma(\{2,3,5\})^2$ gives the following statement
\begin{Proposition}  \label{P235} The nonempty sets $\Polys_{2^b 1^{a}}(\{2,3,5\})$ have size 
as in Table~\ref{V235}.\
\end{Proposition}
\noindent The computation required to carry out Step 3 and thereby 
prove Proposition~\ref{P235} took about two hours.  
\begin{table}[htb]
\[
\begin{array}{r|rrrrrr}
  b & a=0 & a=1 & a=2 & a=3 & a=4 & a=5 \\
  \hline
  0&                            1    &99     &1020   &3100   &3570  &1386 \\
  1&                1927  &18225  &60240  &90640  &64470 &18018\\
  2&                  44967 &227751 &477540 &511200 &279930&64176\\
  3&                 238255&862029 &1347060&1125940&502530&99960\\
  4&                  551944&1567746&1913760&1269160&463470&83034\\
  5&                  745824&1740246&1683180&867600 &246120&40698\\
 6 &                    692476&1364910&1050150&409570 &81690 &12768\\
  7&                  480862&812520 &493440 &146800 &20370 &3360 \\
  8&                 259974&376650 &170850 &38550  &3990  &756  \\
  9&                    112016&138096 &39660  &6020   &420   &84   \\
  10&                39404 &42216  &5520   &380    &      &     \\
   11&                 11520 &11436  &360    &       &      &     \\
   12&                 2751  &2709   &       &       &      &     \\
  13&                  495   &495    &       &       &      &     \\
   14 &                 57    &57     &       &       &      &     \\
   15&                 3     &3      &       &       &      &     \\
\end{array}
\]
\caption{\label{V235} Size of the nonempty sets $\Polys_{2^b1^a}(\{2,3,5\})$.}
\end{table}
The fact that all $a$'s appearing in Table~\ref{V235} 
are at most five is known by Reduction
Bound~\ref{trivbound}, because $\bbP^1(\F_7)$ has only five elements besides $0$, $1$, and $\infty$. 
In contrast, $\bbP^1(\F_{49}) - \{0,1,\infty\}$ has $N(7,2)=47$ elements, corresponding to the bound $2 b + a  \leq 47$.  
Thus our computation identifies many $\Polys_{2^b 1^{a}}(\{2,3,5\})$ as
empty even though the reduction bound  allows them to be non-empty.  

\subsection{Extremal Polynomials}
One of the three elements in $\Polys_{2^{15}1}(\{2,3,5\})$ is 
\begin{eqnarray}
\nonumber s(t) & = &   \left(t^2+6 t+3\right)  \left(3 t^2+6 t+1\right) 
   \left(t^2-6 t+3\right)  \left(3 t^2-6 t+1\right) \cdot \\ 
\nonumber &&    \left(t^2-2 t-5\right)  \left(5 t^2+2 t-1\right)   
 \left(t^2+2 t-5\right)  \left(5 t^2-2 t-1\right)   \cdot \\
\label{big235} &&    \left(t^2-2 t-1\right)  \left(t^2+2 t-1\right) \cdot 
    \left(t^2-6 t-1\right) \left(t^2+6 t-1\right) \cdot \\
 \nonumber &&   \left(3 t^2-2 t-3\right) \left(3 t^2+2 t-3\right) \cdot \left(t^2+1\right) \cdot 
    (t+1).
 \end{eqnarray}
 Its discriminant is $2^{1046} 3^{80} 5^{104}$.  Its
 roots, together with $1$, are visibly invariant under the four-element
 group generated by negation
 and inversion, with minimal invariant factors separated by $\cdot$'s.  
 The other two elements of $\Polys_{2^{15}1}(\{2,3,5\})$ are obtained
 from the given one by  applying the 
 transformations $t \mapsto 1-t$ and $t \mapsto t/(t-1)$.

\section{3-split polynomials}
\label{3polys}
    This section describes how one determines sets $\Polys_{3^c2^b1^{a}}(P)$.
Without loss of generality we restrict to $\cP$ containing $2$ and $3$ throughout this section.  
Assuming that $\Polys_1(P)$ and $\Polys_2(P)$ are known from the 
previous two sections, to complete Step 1 and 2 of the introduction, one
needs to determine $\Polys_3(P)$ and Theorem~\ref{thm3111} gives our method.  We illustrate
the full procedure by determining all $\Polys_{3^c 2^b 1^{a}}({\{2,3\}})$.  

\subsection{Compatible $ABC$ triples and vertices}  Let $T_{3,2,\infty}(\Z^P)$ be the 
set of rational numbers $j = -A/C$ exactly as in Definition~\ref{abc} except that
$A$ and 
$B = -A-C$ are only required to have the respective forms $a x^3$ and $by^2$ with $a,b \in  P^*$.   
For an element $j \in T_{3,2,\infty}(\Z^P)$ the polynomial 
\begin{eqnarray}
\label{jcubic}
S(j,t) & = & 4(j-1) t^3 - 27 j t -27 j
\end{eqnarray}
has discriminant $3^9 j^2/2^4 (j-1)^3$.   Let $c$ be the isomorphism class of the
algebra $\Q[t]/S(j,t)$.    This invariant gives a decomposition
\begin{equation}
\label{3decomp}
T_{3,2,\infty}(\Z^P) = \coprod_c T_{3,2,\infty}(\Z^P)^c.
\end{equation}
This decomposition is used in Theorem~\ref{thm3111} below as one of two aspects of
compatibility.  

To find all $j$ in $T_{3,2,\infty}(\Z^P)$ up to a height bound of $H$, one searches 
as before over $(A,C)$.  Now, however, the search is substantially larger as one
only has $A = a x^3$ with $a \in  P^*$.  
   For our example,
we need the set $T_{3,2,\infty}(\Z^{\{2,3\}})$.  
A three minute search
up to cutoff $10^{11}$ found 81 elements.  Of these, the factorization partition of
\eqref{jcubic} is $3$, $21$, and $1^3$ respectively $54$, $24$, and $3$ times.  
The four $j$'s with \eqref{jcubic} irreducible of largest height come 
from the triples 
\begin{eqnarray*}
 (-73085409, 73085401, 8) &=& (-3^5 67^3, 8466^2 ,2^3), \\
 (128787625, -531440809, 402653184) &=& (505^3,-25053^2,2^{27} 3), \\ 
  (7022735875, -7022744067, 8192) &=& (1915^3,3^1 48383^2, 2^{13}), \\ 
  (67867385039, -67867385042, 3) &=& (4079^3, -2^1 184211^2, 3).
\end{eqnarray*}
All the other elements have height at most 3,501,153.  The completeness
of this $81$-element list dates back to \cite{Co}; it is also 
a subset of the $440$-element set $T_{3,2,\infty}(\Z^{\{2,3,5\}})$ from \cite{Cr} 
cited in the previous section.   The distribution of the $54$ irreducible $j$-invariants according to isomorphism class $c$ 
is quite uneven, and given in Table~\ref{tab3111} below.

\subsection{From the set $T_{3,2,\infty}(\Z^P)$ of $ABC$ triples to the set $\Polys_3(P)$ of degree three vertices in
$\Gamma(P)$}   The current situation is similar to the passage from $T_{\infty,2,\infty}(\Z^P)$ to
$\Polys_2(P)$ but more complicated.    The discriminant of a monic cubic polynomial $s(t) = t^3 + b t^2 + c t + d$
is
\[
\Delta(b,c,d) = -4 b^3 d+b^2 c^2+18 b c d-4 c^3-27 d^2. 
\]
If $s(t)$ is separable, so that $\Delta(b,c,d)$ is nonzero, the 
$j$-invariant is then
\[
j = \frac{4 (b^2-3 c)^3}{27 \Delta(b,c,d)}.  
\]
If one changes $s(t)$ to $m^{-3} s(m t + b)$ the $j$-invariant does not change. 
One can expect $j$-invariants to play a central role in our situation because 
for $j \neq 0,1$, polynomials in $s(t) \in \Q[t]$ with a given $j$-invariant 
are all transforms of each other by fractional linear transformations in
$PGL_2(\Q)$.   

  Let 
\begin{eqnarray*}
F(j,k,y) & = & k \left(j^2 y^3-2 j y^3+3 j y^2-3 j y+1\right)^2-j \left(j y^2-2 y+1\right)^3 \\
             & = & j^2 \left(j^2 k-j^2-4 j k+4 k\right) y^6 + \mbox{(terms of lower order in $y$)}.
\end{eqnarray*}
We say that $\infty$ is a root of $F(j,k,y)$ if the coefficient of $y^6$ is zero.    
This polynomial is important for us because for $j,k \in \Q - \{0,1\}$, 
roots of $F(j,k)$ in $\overline{\Q} \cup \{\infty\}$ are in natural $\Gal(\overline{\Q}/\Q)$-equivariant
bijection with
bijections from roots of $S(j,t)$ to roots of $S(k,t)$.  Note that 
\[
\disc_y(F(j,k,y)) = 2^{22} 3^6 j^{10} (j-1)^{15} k^4 (k-1)^3.
\]
Thus there indeed always six roots when $j,k \in \Q - \{0,1\}$.  
 
\begin{Theorem}   
\label{thm3111} Let $P$ be a finite set of primes containing $2$ and $3$.   Let $c$ be the isomorphism class of a cubic field in $NF_3(P)$ and let
$j \in T_{3,2,\infty}(\Z^P)^c$.    The polynomials in $\MPolys_3(P)^c$ with $j$-invariant $j$  
are among the polynomials
\begin{eqnarray*}
s^{m,n}_{j_0,j_1,j}(t)  & = & \frac{(j-1) (t (n-m)-n)^3+(j-1) j m^3 n^3-j (m n-m t+n t-n)^3}{(m-n)^3}
\end{eqnarray*}
with $j_0$, $j_1 \in (T_{3,2,\infty}(\Z^P))^c \cup \{0\}$.  
Here $m$ and $n$ run over solutions in $\Q \cup \{\infty\}$ of
$F(j,j_0,y)=0$ and $F(j,j_1,y)=0$ respectively.  
\end{Theorem}

\noindent  If $m$ and/or $n$ is $\infty$, one needs to understand the definition of $s^{m,n}_{j_0,j_1,j}(t)$
in a limiting sense.   For example, 
$s_{j_0,j_1,j}^{\infty,n}(t) = t^3 - 3 j n t^2 + 3 j n^2 t  +  j (j-2)  n^3$

\proof

     Let $s(t) = s_\infty(t) = t^3 + b t^2 + c t + d$ be a polynomial in $\MPolys_3(P)^c$.  
 Then one has not only its usual $j$-invariant $j$, but also the
 $j$-invariants $j_0$ and $j_1$ of the transformed polynomials
 \[
 s_0(t) = \frac{t^3}{d} s\left( \frac{1}{t} \right)
 \mbox{ and } s_1(t) = \frac{(t-1)^3}{1+b+c+d}  s \left( \frac{t}{t-1} \right) . 
 \]
 These two new $j$-invariants lie in 
 $(T_{3,2,\infty}(\Z^P))^c \cup \{0\}$, with
 $0$ only possible only if $\disc(s)$ is $-3$ times a square.

    Recovering all possibilities for $s(t)$ from the three invariants $(j_0,j_1,j)$  
is complicated, because thirty-six different polynomials
$x^3 + b x^2 + c x + d \in \overline{\Q}[x]$ give rise to a given
generic $(j_0,j_1,j) \in \overline{\Q}^3$.    
Note that
\[
\disc_t(s^{m,n}_{j_0,j_1,j}(t)) = \frac{2^2 3^3 j^2 (j-1)^2 m^6 n^6}{(m-n)^6}.
\]
For generic $(j_0,j_1,j)$ the thirty-six polynomials are just the thirty-six 
$s_{j_0,j_1,j}^{m,n}(t)$ as $(m,n)$ varies over solutions to $F(j,j_0,m)=F(j,j_1,n)=0$.   
The coordinate relations 
\begin{eqnarray}
m & = &  -\frac{9 \Delta(b,c,d) }{2 \left(b^2-3 c\right) \left(b^2 c+9 b d-6 c^2\right)}, \\
n & = &   -\frac{9 \Delta(b,c,d) }{2 \left(b^2-3 c\right) \left(2 b^3+b^2 c-9 b c+9 b d-6 c^2+27 d\right)}
\end{eqnarray}
let one verify this statement algebraically.  

    There remains the concern that for nongeneric $(j_0,j_1,j)$, there may be cubics
 in $\MPolys_3(\cP)^c$ which are not among the $s_{j_0,j_1,j}^{m,n}(t)$.  
 This indeed happens in the excluded case $j=0$, as discussed just 
 after this proof.  The case $j=1$ is not relevant for the theorem because $S(1,t) = -27 (1+t)$ is not 
an irreducible cubic.   
The cases $m=0$, $n=0$, and $m=n$ arise only when $j=j_0$, $j=j_1$,
and $j_0=j_1$.  Corresponding polynomials in $\MPolys_3(P)^c$ would
have to be stable under $t \mapsto 1/t$, $t \mapsto t/(t-1)$, or $t \mapsto 1-t$ 
respectively.  But there are no such stable polynomials because the commutator
of the possible Galois groups $A_3$ and $S_3$ in $S_3$ does not contain
an element of order two.   Thus, despite the occasional inseparability of
$s^{m,n}_{j_0,j_1,j}(t)$, all polynomials in 
$\MPolys_3(P)^c$ with nonzero $j$-invariant are indeed among the $s^{m,n}_{j_0,j_1,j}(t)$.  
\qed

To get the complete determination of $\MPolys_3(P)$ we need to complement the polynomials
coming directly from Theorem~\ref{thm3111} with the list of polynomials with $j=0$.  A calculation shows 
that there are no 
separable 
polynomials at all with $(j_0,j_1,j)=(0,0,0)$.   So if $j$ is zero, 
at least one of $j_0$ and $j_1$ is nonzero.  So the remaining polynomials 
are in fact just $S_3$-translates of polynomials already found.  

Some further comments clarify Theorem~\ref{thm3111} and how it is used in the
construction of $\MPolys_3(P)$.  Since $j \neq 0$ and the $j_i$
 belong to the
same cubic field, there is a common Galois group, $G \in \{A_3,S_3\}$.    
The polynomials $F(j,j_i,y)$ can factor into irreducibles in three different ways:
\[
F(j,j_i,y) = \left\{ 
\begin{array}{ll}
(\mbox{cubic})(\mbox{quadratic})(\mbox{linear}) & \mbox{($G=S_3$, $j_i \neq 0$),} \\
(\mbox{cubic})(\mbox{linear})(\mbox{linear})(\mbox{linear}) & \mbox{($G=A_3$),} \\
-(s^2-1)^4 (y - \frac{1}{1-s})^3 (y- \frac{1}{1+s})^3 & \mbox{($G=S_3$, $j_i=0$, $j = 1-s^2$).} 
\end{array}
\right.
\]
In the $A_3$ case, always $3^2=9$ of the thirty-six 
$s_{j_0,j_1,j}^{m,n}(t)$ are in $\Q[x]$.  In the $S_3$ case with $w$
zeros among $\{j_0,j_1\}$ there are $2^w$ rational polynomials.
Of course it is trivial to see whether a candidate $s_{j_0,j_1,j}^{m,n}(t)$ from the theorem
is actually in $\MPolys_3(P)$.  Namely, if $m$ and $n$ are both finite then the quantities
\begin{eqnarray*}
s^{m,n}_{j_0,j_1,j_\infty}(0) & = & \frac{n^3 \left(j^2 m^3-2 j m^3+3 j m^2-3 j m+1\right)}{(m-n)^3}, \\
s^{m,n}_{j_0,j_1,j_\infty}(1) & = & \frac{m^3 \left(j^2 n^3-2 j n^3+3 j n^2-3 j n+1\right)}{(m-n)^3}, \\
\disc_t(s_{j_0,j_1,j_\infty}^{m,n}(t)) & = &     \frac{108 (j-1)^3 j^2 m^6
    n^6}{(m-n)^6}
\end{eqnarray*}
need to all be in $\Z^{P \times}$.  When $m$ and/or $n$ is $\infty$, one just uses the limiting
forms of these expressions. 

Table~\ref{tab3111} summarizes the determination of $\Polys_3(\{2,3\})$.
\begin{table}[htb]
\[
\begin{array}{rrl  c r  c |  llrcr}
                   \multicolumn{1}{c}{d} & \multicolumn{1}{c}{D} & \multicolumn{1}{c}{f(t)} & \multicolumn{3}{c|}{T_{3,2,\infty}(\Z^{\{2,3\}})^c}    
                                      &\multicolumn{5}{c}{ |\Polys_3(\{2,3\})^c|} \\
\hline 
                  -6 & -216 & t^3+3 t-2 & \;\;\;\;\;\; &10 && 6(66) &&& = & 396 \\
                  -3 & -972 & t^3-12 & &1 && 6(  &  & 1) & = & 6\\
                  -3 & -972 & t^3-6 & &1 & &6(  &  & 1) & = & 6 \\
                  -3 & -243 & t^3-3 & &6 & &6(13 & + 13 & +4) & = & 180\\
                  -3 & -108 & t^3-2 & &4 &  &6(4 & + 9 &+3) & = & 96  \\
                  -2 & -648 & t^3-3 t-10 && 9 & &6(17) &&& = & 102   \\
                  -1 & -324 & t^3-3 t-4 & &9 &  &6(44) &&& = &264  \\
                  1 & 81 & t^3-3 t-1 & &3 && 6(6 & + 10)& + 2(2)&=&100  \\
                  6 & 1944 & t^3-9 t-6 & & 11 &&  6(58)&&&=& 348\\
                  \cline{5-5} \cline{11-11}
                  &  &   &   &54 &   & &&&&  1498
                 \end{array}
\]
\caption{\label{tab3111} The nine cubic fields $\Q[t]/f(t)$ with discriminant $\pm 2^a 3^b$ and 
associated integers.}
\end{table}
The last block of columns illustrates how a general decomposition of $\Polys_3(P)^c$ into $S_3$-orbits 
appears in the case $P = \{2,3\}$.  
Let $d$ be the square-free integer agreeing with the field discriminant $D$ modulo squares.  If $d \neq 1$
then all orbits have size six.  Orbits are usually indexed by triples of distinct $j$-invariants.  However
for $d = -3$, an unordered triple $(j_0,j_1,0)$ can index up to two orbits and an 
unordered triple $(0,0,j)$ can index up to one orbit.  The contributions from each possibility in the case
$P  = \{2,3\}$ are listed in
order.  For $d=1$, an unordered  triple $(j_0,j_1,j_\infty)$ can index up to 
$9$, $6$, or $1$ orbits, depending on whether it contains $3$, $2$, or $1$ distinct
$j$-invariants.  All orbits again have size six, except for the ones indexed by 
$(j,j,j)$, which have size two.  The contributions from each possibility in the
case $P = \{2,3\}$ are again listed in order.  
 
 \begin{table} 
{\small
 \[
 \begin{array}{rr|rrrrrrrr}
 c & a & b=0 & b=1 & b=2 & b=3 & b=4 & b=5 & b=6 & b=7 \\
 \hline
   0  & 0 &   1     & 169    & 981   & 1723  & 1390 & 630  & 150 & 12 \\
   0    & 1 & 21     & 675    & 2175  & 2559  & 1416 & 486  & 108 & 12 \\
    0   & 2 & 60     & 840    & 1710  & 1200  & 270  &      &     &    \\
   0    & 3 & 40     & 340    & 570   & 340   & 70   &      &     &    \\
 \hline
    1  &  0  & 1498   & 6364   & 10854 & 8788  & 3958 & 1116 & 162 &    \\
    1  & 1 & 4584   & 13632  & 18024 & 11280 & 3600 & 792  & 96  &    \\
    1  & 2 & 4260   & 9900   & 10020 & 4800  & 720  &      &     &    \\
    1  & 3 & 1120   & 2440   & 2040  & 1000  & 160  &      &     &    \\
 \hline
    2  &  0  & 21282  & 37374  & 34008 & 16866 & 4560 & 798  & 72  &    \\
    2  & 1 & 41184  & 62208  & 49872 & 21000 & 3900 & 564  & 48  &    \\
    2  & 2 & 24720  & 33180  & 23160 & 8940  & 900  &      &     &    \\
    2  & 3 & 3960   & 6000   & 3720  & 1680  & 240  &      &     &    \\
 \hline
    3  & 0& 81850  & 95578  & 54942 & 17398 & 2704 & 216  &     &    \\
    3  & 1 & 117288 & 133632 & 71712 & 19800 & 1992 & 120  &     &    \\
    3  & 2 & 49140  & 54660  & 27240 & 7380  & 540  &      &     &    \\
    3  & 3 & 4520   & 6200   & 2760  & 1000  & 160  &      &     &    \\
 \hline
    4  & 0 & 156924 & 144000 & 55692 & 11434 & 1132 & 48   &     &    \\
    4  & 1 & 180822 & 174564 & 64074 & 11004 & 684  & 24   &     &    \\
    4  & 2 & 56910  & 56940  & 19050 & 2760  & 120  &      &     &    \\
    4  & 3 & 3030   & 4020   & 1230  & 220   & 40   &      &     &    \\
 \hline
    5  & 0 & 173110 & 137530 & 38094 & 4848  & 282  &      &     &    \\
    5  & 1 & 167448 & 144552 & 39048 & 3936  & 144  &      &     &    \\
    5  & 2 & 42000  & 37260  & 8880  & 420   &      &      &     &    \\
    5  & 3 & 1240   & 1600   & 360   &       &      &      &     &    \\
 \hline
    6  & 0 & 116552 & 85214  & 18186 & 1392  & 42   &      &     &    \\
    6  & 1 & 95388  & 76440  & 16572 & 1044  & 24   &      &     &    \\
    6  & 2 & 19800  & 15360  & 2820  &       &      &      &     &    \\
    6  & 3 & 560    & 620    & 60    &       &      &      &     &    \\
 \hline
    7  & 0 & 49364  & 33650  & 5622  & 246   &      &      &     &    \\
    7  & 1 & 33576  & 25440  & 4392  & 192   &      &      &     &    \\
    7  & 2 & 5820   & 4140   & 600   &       &      &      &     &    \\
    7  & 3 & 160    & 160    &       &       &      &      &     &    \\
 \hline
    8  & 0  & 12998  & 7916   & 954   & 24    &      &      &     &    \\
    8  & 1 & 6870   & 4914   & 534   & 18    &      &      &     &    \\
    8  & 2 & 960    & 720    & 60    &       &      &      &     &    \\
    8  & 3 & 20     & 20     &       &       &      &      &     &    \\
\hline
    9  & 0 & 1948   & 952    & 54    &       &      &      &     &    \\
    9  & 1 & 648    & 456    &       &       &      &      &     &    \\
    9  & 2 & 60     & 60     &       &       &      &      &     &    \\
\hline
    10 & 0& 162    & 54     &       &       &      &      &     &    \\
    10 & 1 & 24     & 24     &       &       &      &      &     &    \\
\hline
    11 & 0 & 8      & 2      &       &       &      &      &     &    \\
\end{array}
\]
}
\caption{\label{V23} Size of the nonempty sets $\Polys_{3^c2^b1^a}({\{2,3\}})$ }
\end{table}

As an example of the complications associated to $d=-3$, let $c$ be the isomorphism class of 
$\Q[t]/(t^3-12)$.  Then $T_{3,2,\infty}(\Z^{\{2,3\}})^c = \{-24\}$.  
Consider $(j_0,j_1,j) = (0,0,-24).$ 
 Theorem~\ref{thm3111} formally yields four candidates.  
\[
 \begin{array}{c|cc}
 m \setminus n & -1/4 & 1/6  \\
 \hline
 -1/4 & [1] & t^3+6 t^2-6 t+2 \\
 1/6 & t^3-9 t^2+9 t-3 & [1] \\
 \end{array}
 \]
 As always for $(0,0,j)$, only the two candidates coming from $m \neq n$ are
 separable and these are $S_3$-transforms of one another.  In this case,
 both are in $\Polys_3(\{2,3\})^c$. The remaining $S_3$-transforms
 are $2 t^3-3$, $2 t^3-6 t^2+6 t+1$, $3 t^3-2$, and $3 t^3-9 t^2+9 t-1$, 
 accounting for all of $\Polys_3(\{2,3\})^c$.
  
 As an example of complications associated to $d=1$, let $c$ come from $A_3$ field $\Q[t]/(t^3-t-1)$.   
Then $T_{3,2,\infty}(\Z^{\{2,3\}})^c = \{1372/3,4,4/3\}$.  The ordered tuple
 $(j_0,j_1,j) = (1372/3,4,4/3)$ gives nine candidates.  They
 are
 \[
 {\renewcommand{\arraycolsep}{1pt}
 {
 \begin{array}{c|ccc}
 m \setminus n & 1 & 1/2 & \infty \\
 \hline
 3 & 8 t^3-36 t^2+30 t+1 & [125 t^3-225 t^2+75 t+1] & t^3+9 t^2+15 t-1 \\
 3/8 & [125 t^3-300 t^2+180 t-8] & t^3-6 t^2+9 t-1 & 64 t^3-96 t^2+36 t-1 \\
 9/10 & t^3+6 t^2-96 t+8 & 64 t^3-48 t^2-96 t-1 & [125 t^3+75 t^2-120 t+1] 
 \end{array}
 }
 }
 \]
 The three in brackets are rejected and the other six are members of
 $\Polys_3(\{2,3\})^c$.

\subsection{The sets $\Polys_{3^c 2^b 1^a}(\{2,3\})$}
Inductively tabulating cliques in the
graph $\Gamma(\{2,3\})^3$ takes about 15 minutes and yields the following statement.  

\begin{Proposition} \label{Prop23} The nonempty sets $\Polys_{3^c 2^b 1^a}(\{2,3\})$ have size as in
Table~\ref{V23}.
\end{Proposition}
\noindent 
 The largest
 $a$ in Table~\ref{V23} being $3$ agrees with Reduction Bound~\ref{trivbound}, as there are  only $3$ points in
 $\bbP^1(\F_5)$ beyond those already used by $\{0,1,\infty\}$.

\subsection{Extreme polynomials}
One of the two polynomials in $\Polys_{3^{11} 1}(\{2,3\})$ has already been given 
in \eqref{big23}.  It is stable under 
 the group $A_3$ of even permutations of the cusps $\{0,1,\infty\}$.  
 The minimal $A_3$-stable factors are given between adjacent $\cdot$'s in 
 \eqref{big23}.

 \section{$f$-split polynomials with $f \geq 4$}
 \label{genpolys}
     In the previous two sections we have used $w$-invariants in $T_{\infty,2,\infty}(\Z^P)$
      to construct sets $\Polys_2(P)$ 
and $j$-invariants in $T_{3,2,\infty}(\Z^P)$ 
to construct sets $\Polys_3(P)$.    For degrees $d \geq 4$, we follow the three-step approach
of \S\ref{33} to determining $\Polys_d(P)$.  

\subsection{Excellent $P$-units and $\Polys_d(P)$} 

    Let $K_1$,\dots, $K_m$ be a list of degree $d$ number fields
  unramified outside $P$ such that every isomorphism class of such 
  fields appears once.      Let $U_i$ be the $P$-unit group of $K_i$.  
  The finitely generated group $U_i$ is well-understood and there
  are algorithms to produce generators.  A $P$-unit $u$ is called 
  an {\em exceptional $P$-unit} if $1-u$ is also a unit.   Exceptional $P$-units have
  been the subject of many studies, e.g.\ \cite{Ni}.  
  
      Let $s_u(t)$ be the characteristic polynomial of
   a $P$-unit $u$.  So $s_u(t)$ is a monic polynomial of degree $d$ 
   in $\Z^P[t]$ with constant term $s_u(0)$ in $\Z^{P \times}$.  The unit $u$  is exceptional if
   and only if also $s_u(1) \in \Z^{P \times}$.   We say it is an {\em excellent $P$-unit} if
   furthermore $\disc(s_u(t)) \in \Z^{P \times}$.   All elements of $\MPolys_d(P)$ arise in
   this way as characteristic polynomials of excellent units.   
   
       As the example of this section, we take $P=\{2\}$.  The relevant set 
  $NF_d(\{2\})$ of number fields is known in degrees $d \leq 15$ \cite{Jo,Database}.    In fact, for $d=1$, $2$, $4$, $8$
  one has $|NF_d(\{2\})|=1$, $3$, $7$, $36$ and otherwise $NF_d(\{2\}) = 0$.
  The fields in question are all totally ramified
  at $2$.  This implies that an $S_3$-orbit of polynomials in $\Polys_d(\{2\})$ takes
  one of the following forms.  First, the orbit may contain just three polynomials, one
  of which is palindromic.  In this case the palindromic polynomial $s(t)$ is the
  unique member of the orbit with $2|s(1)$.  
 Second, the orbit may contain six polynomials, none palindromic.  In this case,
  exactly two of the polynomials $s_i(t)$ satisfy $2|s_i(1)$.   They 
  are related by $t^d s_1(1/t) = \pm s_2(t)$.

\begin{table}[htb]
{\small
{\renewcommand{\arraycolsep}{3pt}
\[
\begin{array}{| lrlllcrrccc |}
\hline
   G & d & \lambda_\infty & \lambda_3 & s(t) &\mbox{Pal}&s(1)&c& N_1 & N_2 & N_4 \\
   \hline
   \multicolumn{10}{c}{\;}  \\
   \hline
   C_1&1&1&1&                 t+1&\bullet &2 & 1 &0  & 7 & 36 \\
                    \hline
                    \multicolumn{10}{c}{\;} \\
                    \hline
   C_2  &8&11&2&                   t^2+6 t+1&\bullet &8 & 2 & 1 &0  & 7 \\
   &&&&                     t^2-6 t+1 &\bullet &-4 & 2 &  1 & 0 & 7 \\
       &&&&                 t^2-2 t-1&&-2 & 1 & 2 & 2 & 13 \\
                    \hline
      C_2 &-4&2&2&                 t^2+1&\bullet &2 & 1 & 1 & 2 & 14 \\
                    \hline
       C_2&-8&2&11&&&&&&& \\
                    \hline
                    \multicolumn{10}{c}{\;}  \\
                    \hline
      V & 256 & 2 & 22&                  t^4+1&\bullet &2 & 1 &  1 & 1 & 4 \\
       &&&&                 t^4+6 t^2+1& &8 & 8 & 1 & 1 & 2 \\
                    \hline
       C_4 & 2048 & 1111 & 4        &        t^4-4 t^3-26 t^2-4 t+1&\bullet &-32 & 64 & 1 & 2 & 2 \\
              &&&&             t^4+4 t^3-26 t^2+4 t+1 &\bullet &-16 & 64 & 1 & 2 & 2 \\
       &&&&             t^4+28 t^3+70 t^2+28 t+1 &\bullet& 128 & 512 & 1 & 1 & 0 \\
       &&&&             t^4-28 t^3+70 t^2-28 t+1 &\bullet& 16 & 512 & 1 & 1 & 0 \\
       &&&&             t^4-4 t^3-6 t^2+4 t+1 & & -4 & 8 & 1 & 3 & 3 \\
       &&&&             t^4-4 t^3-2 t^2+12 t+1& &8 & 8 &  2 & 1 & 4 \\
              &&&&               t^4+4 t^3-2 t^2-12 t+1 &&8 & 8 & 1 & 1 & 4 \\
       &&&&             t^4-4 t^3-2 t^2+4 t-1 && -2 & 1 &  2 & 2 & 7 \\
           &&&&              t^4+4 t^3-2 t^2-4 t-1& &-2 & 1 & 1 & 1 & 5 \\
       &&&&             t^4-148 t^3+102 t^2-20 t+1& &-64 & 512 &  0&0  & 1 \\
       &&&&             t^4-20 t^3+34 t^2-12 t+1 &&4 & 32 &  0 & 2 &  0\\
                                          \hline
                                       C_4 & 2048 & 22 & 4 && &&&&& \\
                    \hline
                    D_4 & -2048 & 211 & 22 &
              t^4+12 t^3+6 t^2+12 t+1 &\bullet&32 & 64 & 1 & 2 & 2 \\
      &&&&               t^4-12 t^3+6 t^2-12 t+1 &\bullet& -16 & 64 & 1 & 2 & 2 \\
     &&&&               t^4-4 t^3+6 t^2-4 t-1 &&-2 & 1 & 1 & 2 & 5 \\
                    \hline
                    D_4 & -2048 & 211 & 4 &
                             t^4-4 t^3-2 t^2-4 t+1 &\bullet&-8 & 8 & 1 & 2 & 6 \\
                                   &&&&              t^4+4 t^3-2 t^2+4 t+1&\bullet & 8 & 8 & 1 & 2 & 6 \\
      &&&&              t^4+20 t^3-26 t^2+20 t+1&\bullet &16 & 512 & 1 & 1 & 0  \\
       &&&&             t^4-20 t^3-26 t^2-20 t+1&\bullet &-64 & 512 & 1 & 1 & 2 \\
      &&&&              t^4-2 t^2-1& &-2 & 1 & 1 & 2 & 8 \\
      &&&&             t^4-12 t^3+10 t^2-4 t+1 & &-4 & 8 & 1 & 2 & 2 \\
                \hline
                D_4 & 2048 & 22 & 211 &&&&&&& \\
                    \hline
        
    D_4   &512 & 22 & 4 &             t^4-4 t^3+22 t^2-4 t+1&\bullet &16 & 64 & 1 & 1 & 2 \\
       &&&&             t^4+4 t^3+22 t^2+4 t+1&\bullet &32 & 64 &  1 & 1 & 2 \\
       &&&&             t^4+4 t^2-4 t+1 &&2 & 1 & 1 & 1 & 4 \\
                    \hline
 \end{array}
\]
}
}
\caption{\label{Irred2polys} Information on the sets $\Polys_d(\{2\})$ for $d \in \{1,2,4\}$  }
\end{table}

   Table~\ref{Irred2polys} describes the sets $\Polys_d({\{2\}})$ for $d \in \{1,2,4\}$
   by listing polynomials representing $S_3$-orbits.   
   The fields defined by $t^2+2$,  $t^4+4 t^2+2$, and $t^4+2$ yield no polynomials 
   at all.   Their $3$-adic factorization partitions $\lambda_3$ are respectively $1^2$, $4$, and $21^2$.  
   Thus, in the case of $t^2+2$ and $t^4+2$, the nonexistence of polynomials
   follows from Reduction Bound~\ref{trivbound}.     The column $\lambda_\infty$ gives the
   splitting over $\R$.   The rank of $U_i$ is the number of parts of $\lambda_\infty$,
   and, as expected, more parts are correlated with more polynomials.  Palindromic
   polynomials contributing
   three  and nonpalindromic polynomials contributing six, one gets
   $|\Polys_1(\{2\})| = 3$,   $|\Polys_2(\{2\})| = 15$, and $|\Polys_4(\{2\})| \geq 108$. 
   We have taken the computation far enough that it seems unlikely
   that $\Polys_4(\{2\})$ contains polynomials beyond those we 
   have found.    The polynomial discriminants are $Dc^2$, and 
   the largest magnitude $2^{29}$ arises in five orbits.

\subsection{The sets $\Polys_{4^d 2^b 1^a}(\{2\})$}  

          The last three columns of Table~\ref{Irred2polys} indicate the nature of the known part of the graph
   $\Gamma(\{2\})^4=\Gamma(\{2\})^7$.   For each polynomial, 
   the number of neighbors of a given degree $d$ is given as $N_d$.  
   As one would expect in general, the number of neighbors tends to
   decrease as the largest coefficient of the polynomial increases.   
   Carrying out Step~3 as in  \S\ref{twostep} takes less than a second and 
   yields the following result.  
   
   \begin{Proposition} 
   The sets $\Polys_{4^d 2^b 1^a}(\{2\})$ are at least as
   large as indicated
   in Table~\ref{V2}, with equality if $|\Polys_{4}(\{2\})| = 108$.
   \end{Proposition}

\begin{table}
\[
 \begin{array}{l|rrrrr}
         \multicolumn{5}{c}{|\Polys_{4^d 2^b}(\{2\})|}  \\
          b  &        d: \; 0 & 1 & 2 & 3 & 4  \\
                   \hline
           0&   1 & 108 & 177 & 144 & 42 \\
           1&        15 & 162 & 93 & 30 &    \\
           2&        9 & 30 & 21 & 6 &   \\
           3&        3 & 6 & 3 &  &    \\
              \end{array}
\;\;\;\;\;\;\;\;\;\;\;\;
   \begin{array}{l|rrrrr}
                 \multicolumn{5}{c}{|\Polys_{4^d 2^b 1}(\{2\})|}  \\
                  b  &        d: \; 0 & 1 & 2 & 3 & 4  \\
                  \hline
           0&      3 & 108 & 129 & 90 & 24  \\
           1&        21 & 156 & 63 & 18 &   \\
           2&        9 & 18 & 9 &  & \\
          3&         3 & 6 & 3 &  &   
                  \end{array}
\]
 \caption{\label{V2} Size of the nonempty sets $\Polys_{4^d2^b1^a}(\{2\})$,
 assuming  $|\Polys_4(\{2\})| = 108$. }
\end{table}

\subsection{Extremal polynomials}
As an extreme example, the palindromic polynomial
\begin{eqnarray*}
s(t) & = & (t+1) \left(t^2+1\right) \left(t^2-2 t-1\right) \left(t^2+2 t-1\right) \\
&& \qquad    \left(t^4-4 t^3-6 t^2+4 t+1\right) \left(t^4+4 t^3-6 t^2-4 t+1\right)
\end{eqnarray*}
has discriminant is $-2^{184}$.  Its $S_3$ orbit in $\Polys_{4^2 2^3 1}(\{2\})$
corresponds to the bottom right $3$ in Table~\ref{V2}.

\section{Large degree polynomials}
\label{largedegree}
   In this section, we enter the third regime of \S\ref{33}: the systematic construction
of polynomials in $\Polys_\kappa(P)$ in settings where complete determination
of $\Polys_\kappa(P)$ is
well out of reach.   Each subsection focuses on degree $k$ polynomials, 
without pursuing details about their factorization, thus on the sets
$\Polys_{[k]}(P) = \coprod_{\kappa \vdash k} \Polys_\kappa(P)$.
\subsection{Cyclotomic polynomials}
\label{cyclo}
    The following simple result supports the main conjecture of
 \cite{RV}.    

\begin{Proposition} \label{cycloprop} Let $P$ be a finite set of primes containing $2$
and at least one odd prime.  Let $\Polys_{[k]}^{\rm cyclo}(P)$
be the subset of $\Polys_{[k]}(P)$ consisting of products of 
cyclotomic polynomials.   Then $\lim_{k \rightarrow \infty} |\Polys_{[k]}^{\rm cyclo}(\cP)| = \infty$.  
\end{Proposition}

\proof  Let $P^\star$ denotes the set of all 
integers greater than one which are divisible only by primes
of $P$.    For $i \in P^\star$, let $\Phi_i(t)$ be the corresponding
cyclotomic polynomial, of degree $\phi(i)$.  
Then 
\begin{equation}
\label{genfunct}
\sum_{k=0}^\infty |\Polys_{[k]}^{\rm cyclo}(P)| x^k = \prod_{i \in P^{\star}} (1 + x^{\phi(i)}). 
\end{equation}
To treat the sets $P$ appearing in the proposition, we first consider
the case $P = \{2\}$.  One has $\{2\}^\star = \{2,4,\dots,2^j,\dots\}$ and
$\phi(2^j) = 2^{j-1}$.  Expanding the product $(1+x)(1+x^2)(1+x^4) \cdots$,
Equation~\eqref{genfunct} becomes
\begin{equation}
\label{genfunct2}
\sum_{k=0}^\infty  |\Polys_{[k]}^{\rm cyclo}(\{2\})| x^k = 1 + x + x^2+ x^3 + \cdots  
\end{equation}
For $P$ as in the theorem, $\sum |\Polys_{[k]}^{\rm cyclo}(P)| x^k$ is
an infinite sum of $x^j f(x)$ with 
$f(x)$ as in \eqref{genfunct2}.    Thus 
in fact $|\Polys_{[k]}^{\rm cyclo}(\cP)|$ grows monotonically
to $\infty$.  \qed

A numerical example of particular relevance to Hurwitz number fields is 
$P = \{2,3,5\}$.  Then
\begin{eqnarray}
\nonumber \lefteqn{\sum_{\rho=0}^\infty |\Polys^{\rm cyclo}_{[k]}(\{2,3,5\})| x^k } \\
\label{genfunct235}  & = & (1+x)(1+x^2)^3(1+x^4)^4(1+x^6)^2 \cdots \\
\nonumber &=& 1 + x + 3 x^2 + 3 x^3 + 7 x^4 + \cdots + 
3361607445659519 x^{1000} + \cdots 
\end{eqnarray}
We will return to this generating function in \S\ref{largerdegrees}.  

\subsection{Fractal polynomials}   
\label{fract}
A recursive three-point cover is a rational function $F(t) \in \Q(t)$ with all critical values in 
$\{0,1,\infty\}$ and $F(\{0,1,\infty\}) \subseteq \{0,1,\infty\}$.    It has bad reduction
within $P$ if one can write $F(t) = u f(t)/g(t)$ with $f(t)$ and $g(t)$ compatible
polynomials in $\Polys(P)$ and $u \in \Z^{P \times}$.    
   Recursive three-point covers with bad reduction 
within $P$ are closed under composition.  

    The degree $1$ recursive three-point covers form the group $S_3 = \langle 1-t, 1/t \rangle$ 
and have bad reduction set $P = \{\}$.  
Other simple examples are  $F(t) = t^p$ for 
a prime $p$ with bad reduction set $\{p\}$.  
Combining just these via composition one already has a large collection
of recursive three-point covers with solvable monodromy group \cite{RFract}.   One can 
easily extract many other recursive three-point
covers from the literature.  As an example
coming from trinomials, $F_m(t) = t^m/(m t + 1 - m)$ 
has monodromy group $S_m$ and bad reduction 
exactly at the primes dividing $m(m-1)$.   
From the definitions, one has the following fact:

\begin{PullbackConstruction} \label{fractprop} Let $F(t) = u f(t)/g(t)$ be a recursive three-point cover of 
degree $m$ and bad reduction within $P$.  Let $s(t) \in \Polys_{[k]}(P)$.    Then
the pullback $s(F(t)) g(t)^k$ is a scalar multiple of a polynomial in $\Polys_{[mk]}(P)$.   
\end{PullbackConstruction}

\begin{figure}[htb]
\includegraphics[width=4in]{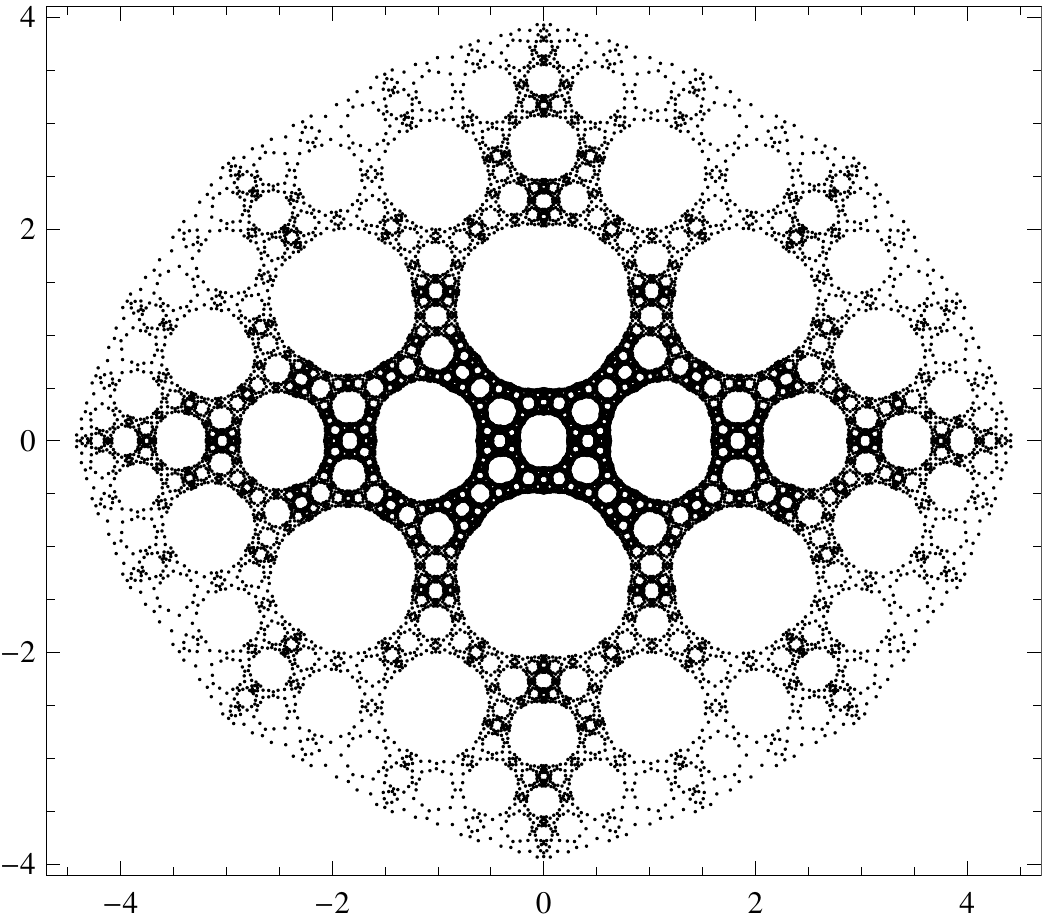}
\caption{\label{fractalroots} The set $R_{8,1}$ consisting of the $2^{15} = 32768$ roots of the specialization polynomial $s_{8,1}(t)$ with bad reduction only at $2$. }
\end{figure}

\noindent We use the word ``fractal'' because when one constructs polynomials
by iterative pullback, their roots tend to have a fractal appearance, as in Figure~\ref{fractalroots}.
 
To explain the source of Figure~\ref{fractalroots}, and also as an example of using the pullback construction iteratively, we prove
 the following complement to Proposition~\ref{cycloprop}.
 \begin{Proposition} \label{fractprop2} The sets $\Polys_{[k]}(\{2\})$ can be arbitrarily large.
 \end{Proposition}
\proof Consider quartic recursive three-point cover
 \[
 F(t) = \frac{-(t-1)^2 (t+1)^2}{4t^2} = \frac{-(t^2+1)^2}{4t^2} + 1 = -\frac{\left(t^2-2 t-1\right) \left(t^2+2 t-1\right)}{4 t^2} -1.
  \]
  Its bad reduction set is $\{2\}$.  
 Some preimages are as follows, with Galois orbits separated by semi-colons:
 \begin{align*}
 F^{-1}(0) & = \{-1; \; 1\}, &   F^{-1}(1) & = \{i,-i\}, \\
 F^{-1}(\infty) & = \{0; \; \infty\}, & F^{-1}(-1) & =  \{-1-\sqrt{2},-1+\sqrt{2}; \; 1-\sqrt{2}, 1 + \sqrt{2} \}.
 \end{align*}
 Let 
 \begin{align*}
 R_{1,-1} & = \{-1 \pm \sqrt{2} \}, &
 R_{1,0} & = \{\pm i\}, &
 R_{1,1} & = \{1 \pm \sqrt{2} \}, &
 \mbox{and } R_{i,j} & = F^{1-i}(R_{1,j}).
 \end{align*}
   The situation is summarized by the following
 diagrammatic description of the action of $F$ on the entire iterated preimage 
 of $\infty$:
 \[
 \begin{array}{ccccc}
 \vdots && \vdots &&  \vdots \\
 \downarrow && \downarrow &&  \downarrow \\
 R_{2,1} && R_{2,-1} &&  R_{2,0} \\
  \downarrow && \downarrow &&  \downarrow \\
 \!\!\! \{1 \pm \sqrt{2}\} \!\!\! && \!\!\! \{-1 \pm \sqrt{2}\} \!\!\! &&  \{\pm i\} \\
  & \!\!\! \searrow \;\;  \swarrow \!\!\! &&& \downarrow \\
  & \{-1\} && & \{1\}    \\
  && \!\!\!\!\!\!\!\!\!\searrow & \swarrow \!\!&\\
  &&  \multicolumn{2}{c}{\{0\}} &\\
  &&  \multicolumn{2}{c}{\downarrow} &\\
  &&  \multicolumn{2}{c}{\{\infty\}} &\\
  &&  \multicolumn{2}{c}{\circlearrowleft} & 
 \end{array}
\]
Note that the critical values $0$, $1$, and $\infty$ have two preimages each while all other values have four preimages.

 For  $i \in \Z_{\geq 1}$ and $j \in \{-1,0,1\}$, let $s_{i,j}(t) \in \Polys_{[2^{2i-1}]}(\{2\})$ be the polynomial  
 with roots $R_{i,j}$.    Products
  of the form $s_{1,j_1}(t) \cdots s_{w,j_w}(t)$ give $3^{w}$ distinct 
  polynomials in $\Polys_{[k]}({2})$  of the same degree $k = 2(4^w-1)/3$.  \qed

 \section{Specialization sets and number field construction}  
 \label{specialization} In this section, we sketch how the sets $\Polys_\kappa(\cP)$ are useful in constructing interesting number fields,
 focusing on two representative families of examples. 
 
 \subsection{Sets $U_\nu(R)$}
   Let $\nu = (\nu_1,\dots,\nu_r)$ be a list of positive integers.  In this section, we assume 
$\nu_{r-2}=\nu_{r-1}=\nu_r=1$ and these indices play a completely passive role.   In the next section, we remove 
this assumption and the last three indices then take on an active role on the same footing with the other
indices.  Without loss of generality, we generally focus on the case where the $\nu_i$ are weakly decreasing,
and use abbreviations such as $21^3 = (2,1,1,1)$.   
  
For any commutative ring $R$, define
$U_\nu(R)$ to be the set of tuples $(s_1(t),\dots,s_{r-3}(t))$ where $s_i(t)$ is a monic degree $\nu_i$ polynomial
in $R[t]$ and the discriminant of 
\[
s_1(t) \cdots s_{r-3}(t) t (t-1)
\]
 is in the group of invertible elements $R^\times$.    To be more explicit,
write $k = \sum_{i=1}^{r-3} \nu_i$ and  
\[
s_i(t) = t^{\nu_i} + u_{i,1} t^{\nu_i-1} + \cdots + u_{i,\nu_i-1} t + u_{i,\nu_i}.
\]
Then the lexicographically-ordered coordinates $u_{1,1}$, \dots, $u_{r-3,\nu_{r-3}}$ realize
$U_\nu(R)$ as a subset of $R^k$.     

The sets $U_\nu(\Z^\cP)$ can be built in 
a straightforward fashion from the sets $\MPolys_\kappa(\cP)$ with
$\kappa$ running over refinements of the partition $(\nu_1,\dots,\nu_{r-3})$.
For example, $U_{1^{k+3}}(\Z^P)$ consists of tuples $(s_1(t),\dots,s_k(t))$
having product $s_1(t) \cdots s_k(t)$ in $\MPolys_{1^k}(\cP)$.    It is thus 
trivially built from $\MPolys_{1^k}(\cP)$, but $k!$ times as big.    As
another example, $U_{k1^3}(\Z^P) = \MPolys_{[k]}(P) = \coprod_{\kappa \vdash k} \MPolys_\kappa(P)$.  In general,
the construction of $U_\nu(\Z^P)$ from $\MPolys_\kappa(P)$ is similar, but
combinatorially more complicated than the two simple extreme
cases just presented.    

\subsection{The scheme $U_\nu$} The object $U_\nu$ itself is an affine scheme, smooth
and of relative dimension $k$ over $\Spec(\Z)$.    We have a focused in Sections~\ref{intro}-\ref{largedegree}
on the sets $\Polys_\kappa(P)$ because of their relatively small size and their direct connection
to graph theory.  However the close variants $U_\nu(\Z^P)$ should be understood 
as the sets of true interest in the application.

\begin{figure}[htb]
\begin{center}
\includegraphics[width=4.8in]{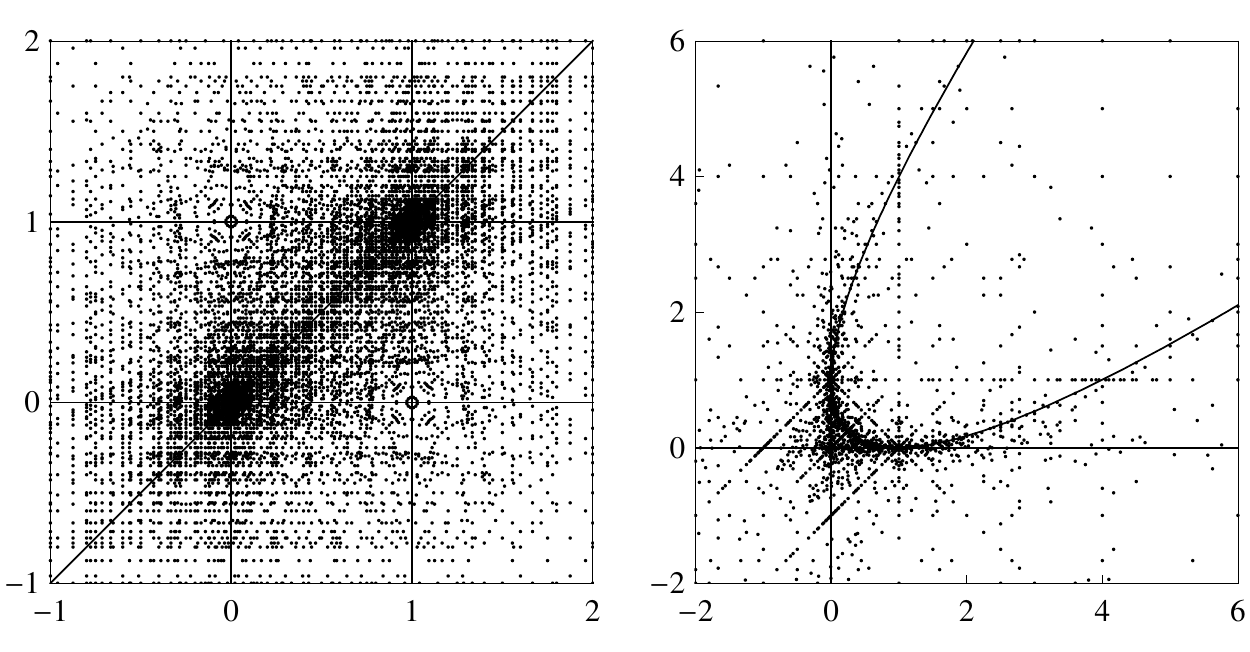}
\end{center}
\caption{\label{pictboth}  $U_{1^5}(\Z^{\{2,3,5,7\}})  \subset U_{1^5}(\R)$ and  $U_{21^3}(\Z^{\{2,3,5\}})  \subset   U_{21^3}(\R)$.}
\end{figure}

The sets $U_\nu(\Z^P)$ fit into standard geometrical considerations much better than
the $\Polys_\kappa(P)$ do.  For 
example, $U_\nu(\Z^P)$ lies in the $k$-dimensional real manifold $U_\nu(\R)$,
while similar oversets are not as natural for $\Polys_\kappa(P)$.   Figure~\ref{pictboth} 
draws examples, directly related to Sections~\ref{1polys} and \ref{2polys}.
In each case, $U_\nu(\R)$ is the complement in $\R^2$ of the drawn curves.

\subsection{Covers} 
The fundamental groups of the complex manifolds 
$U_\nu(\C)$ are braid groups.  
Katz's theory \cite{Ka}
of rigid local systems gives a whole hierarchy of 
covers of $U_{\nu}$ \cite{RABC,RM05}.
The theory of Hurwitz varieties as presented in  \cite{BR} likewise gives 
a another whole hierarchy of  covers \cite{RV,RHNF}.   
In both cases, the covers have a topological
description over $\C$, and this description
forms the starting point of a more arithmetic
description over $\Z$.   In each case,
the datum defining a cover determines
a finite set $W$ of primes.  The cover
is then unramified except in 
characteristics $p \in W$.  
For these bad characteristics, the cover
is typically wildly ramified.   

Rather
than enter theoretically into these two theories,
we discuss next two representative examples,
both with $\nu = 21^3$ for uniformity.  
As coordinates, we work with $(u_{1,1},u_{1,2}) = ((1-v-u)/u,v/u)$
so that the right half of Figure~\ref{pictboth} is the $u$-$v$ plane.  
Our current coordinates are related to the quantities of
Section~\ref{2polys} by $u = u_0/u_1$ and $v = u_\infty/u_1$.  

\subsection{A Katz cover}
\label{Katz}
The last half of \cite{RM05} considers two Katz covers
with bad reduction set $W = \{2,3\}$.  
The smaller of the two is
captured by the explicit polynomial 
\begin{eqnarray*}
\lefteqn{f_{27}(u,v,x) =} \\
&& \left(x^3 - 3 d x + 2 e d\right) \left(x^6    -15 d x^4  +40 e d x^3 -45
    d^2 x^2 +24 e d^2 x  -32 e^2 d^2+27 d^3 \right)^4 \\ && -432 v d \left(x^4 - 6 d x^2 + 8 e d x-3 d^2\right)^6,
\end{eqnarray*}
with abbreviations $d=u^2 + v^2 - 2 u v - 2 u - 2 v + 1$ and $e = u+v-1$.  
The polynomial discriminant factors, 
\[
D_{27}(u,v) = 2^{840} 3^{270} u^{102} v^{126} d^{234}.  
\]
The Galois group of $f_{27}(u,v,x) \in \Z[u,v][x]$ is
the orthogonal group $O_6^-(\F_2) \subset S_{27}$ 
of order $51,840 = 2^7 \, 3^4 \, 5$.  The specialization
set $U_{21^3}(\Z^{\{2,3\}})$ has order $60 + 169 = 229$
from Table~\ref{V23}.   This specialization process
produces 193 number fields with
Galois group $O_6^-(\F_2)$, 15 with Galois group 
the index two simple group $O_6^-(\F_2)^+$, 
and other number fields with 
various smaller Galois groups \cite{RM05}.     
 
Covers in the Katz hierarchy typically
yield Lie-type Galois groups, like in 
this example, with bad reduction
set $W$ containing at least two primes.   By 
varying the Katz cover, a single fixed 
specialization point $u \in  U_\nu(\Z^W)$
with $|W| \geq 2$
can be expected to yield infinitely many 
different fields ramified within $W$.

\subsection{A Hurwitz cover}
\label{Hurwitz}
Many Hurwitz covers of $U_{21^3}$
with bad reduction set $W = \{2,3,5\}$
are studied in \cite{RHNF}.   One such cover
 has degree $36$ and can be given via equations
as follows.    The cover $X$ can be given 
coordinates $x$ and $y$ so that the map to 
$U_{21^3}$ takes the form   
\begin{eqnarray*}
u&=& -\frac{32 \left(x^2-2 y\right)^5}{27 y (x-y)^4 \left(8 x^3-x^2-18 x
    y+27 y^2+2 y\right)}, \\
 v   &=& -\frac{\left(4 x^3-x^2+18 x y+27 y^2-4 y\right)^2
    \left(2 x^4-5 x^2 y+6 x y^2-y^3+2 y^2\right)}{27 y (x-y)^4 \left(8
    x^3-x^2-18 x y+27 y^2+2 y\right)}.
\end{eqnarray*}
Eliminating $y$ gives $f_{36}(u,v,x) \in \Z[u,v][x]$ with $x$-degree
$36$ and $1125$ terms.  Its discriminant is 
\[
D_{36}(u,v) = -2^{337} 3^{513} 5^{750} u^{53} v^{13} d^{22} C(u,v)^2,
\]
with the complicated polynomial $C(u,v) \in \Z[u,v]$ not contributing 
to field discriminants of specializations.    The specialization
set $U_{21^3}(\Z^{\{2,3,5\}})$, drawn as the right half of Figure~\ref{pictboth},  has order $1927 + 1020 = 2947$
from Table~\ref{V235}.   The specialization process
produces 2652 number fields with
Galois group $S_{36}$, 42 number fields with Galois group
$A_{36}$, and others with various smaller Galois groups,
all with bad reduction set exactly $\{2,3,5\}$.

Covers in the Hurwitz hierarchy typically
yield alternating or symmetric groups, like in 
this example, with bad reduction
set $W$ containing all primes dividing the
order of some nonabelian finite simple group, thus
at least three primes.  Here again, 
by varying the cover, 
a single fixed specialization point $u \in  U_\nu(\Z^W)$
can give many different fields ramified within
$W$.    

\subsection{Constraining wild ramification}  Let $K_u$ 
be an algebra obtained by specializing a cover with
bad reduction set $W$ at a point $u$ with
bad reduction set $P$.  Then the typical
behavior of $p$-adic ramification in $K_u$ is as follows:
\[
\begin{array}{l | ll}
  & p \in P & p \not \in P \\
\hline
p \in W & \mbox{very wild} & \mbox{slightly wild} \\
p \not \in W & \mbox{tame} & \mbox{none} \\
\end{array}
\]
To illustrate the distinction between ``very wild'' and ``slightly wild'', 
we specialize the 
Katz cover $f_{27}(u,v,x)$ and the Hurwitz cover
$f_{36}(u,v,x)$ at the $15$-element set 
$U_{21^3}(\Z^{\{2\}})$ appearing as black
vertices in Figure~\ref{twograph}:

{
{\renewcommand{\arraycolsep}{3pt}
\[
\begin{array}{rrr|rrr|rrr}
u_0 & u_1 & u_\infty &  \multicolumn{3}{c|}{d_{27}(u,v)} &  \multicolumn{3}{c}{d_{36}(u,v)}  \\
\hline
8&1&1& 2^{ 88} 3^{32  } &   2^{98 } 3^{ 36 } & \bullet 2^{86 } 3^{34  } & -2^{127   } 3^{ 39   } 5^{30   } &  -2^{127   } 3^{ 43   } 5^{30   } &  -2^{ 118  } 3^{ 39   } 5^{ 30  }  \\
-4&1&1& 2^{80 } 3^{ 36 } &   2^{84 } 3^{36  } &   2^{ 80} 3^{ 36 } &2^{ 114  } 3^{ 39   } 5^{ 30  } &  2^{ 112  } 3^{ 39   } 5^{ 30  } & \bullet{-2^{ 100  } 3^{ 25   } 5^{  28 }}  \\
2&1&1& 2^{88 } 3^{ 36 } &   2^{ 88} 3^{ 30 } &   2^{82 } 3^{ 30 } &-2^{ 115  } 3^{ 39   } 5^{  30 } &  -2^{121   } 3^{ 39   } 5^{ 30  } &  -2^{118   } 3^{39    } 5^{30   }  \\
-2&-1&1& 2^{98 } 3^{ 36 } &   2^{102 } 3^{ 34 } &   2^{96 } 3^{ 32 } &-2^{137   } 3^{ 27   } 5^{ 30  } &  2^{ 137  } 3^{  39  } 5^{ 30  } &  2^{134   } 3^{ 39   } 5^{ 30  }  \\
2&-1&1& 2^{98 } 3^{ 36 } &   2^{102 } 3^{36  } &   2^{96 } 3^{ 36 } & 2^{137   } 3^{ 35   } 5^{ 30  } &  -2^{ 137  } 3^{ 27   } 5^{ 30  } &  2^{ 134  } 3^{  39  } 5^{ 30  } 
\end{array}
\]
}
}

\noindent   
On
a given row starting with $(u_0,u_1,u_\infty)$, the $(u,v)$ in the 
second and third blocks are, in order, $(u_0,u_\infty)/u_1$, $(u_1,u_0)/u_\infty$, and 
$(u_\infty,u_1)/u_0$.

Field discriminants of the specializations are as indicated by the table.   
 As $(u,v)$ runs over all of $U_{21^3}(\Z^{\{2,3\}})$ one gets discriminants
$d_{27}(u,v) = 2^a 3^b$ with $a \in [24,102]$ and $b \in [26,60]$.    Restricting to $U_{21^3}(\Z^{\{2\}})$,
the maximum $a$ appearing is not reduced at all, while the maximum $b$ is reduced from $60$ to $36$.  
Similarly, as  $(u,v)$ runs over all of $U_{21^3}(\Z^{\{2,3,5\}})$, one gets discriminants
$d_{36}(u,v) = \pm 2^a 3^b 5^c$ with $a \in [40,137]$, $b \in [7,72]$, and $c \in [18,61]$.   
Restricting to  $U_{21^3}(\Z^{\{2\}})$,
$a_{\rm max}$  is not reduced at all, while $b_{\rm max}$ is reduced from $72$ to $39$
and $c_{\rm max}$ is reduced from $61$ to $30$.    This distinction between ``very wild'' and ``slightly wild''
makes all the sets $\Polys_\kappa(P)$ of interest in the applications, not just the 
ones where $P$ is large enough to contain the bad reduction set
 $W$ of a cover.  

\subsection{Explicit examples} To give completely explicit examples of number fields constructed using the 
specialization points, we continue the previous
subsection.   The fifteen specializations of $f_{27}(u,v,x)$ in the table all have Galois group 
$O_{6}^-(\F_2)$ except the bulleted one, which has Galois group the simple index
two subgroup.   A presentation for this field $K_{27,1/8,1/8}$ is $\Q[x]/g_{27}(x)$ with  
\begin{eqnarray*}
g_{27}(x) & = & x^{27}-9 x^{26}+21 x^{25}+53 x^{24}-288 x^{23}+1628 x^{21}-1164 x^{20}-5409
    x^{19} \\ && +5681 x^{18}+12159 x^{17}-14793 x^{16}-20548 x^{15}+25764
    x^{14}+30324 x^{13} \\  
    && -36220 x^{12}-42249 x^{11}+48465 x^{10}+50819 x^9-61773
    x^8-44220 x^7 \\ && 
    +64172 x^6+23712 x^5-48024 x^4-5725 x^3+22509 x^2+147 x-5045. 
\end{eqnarray*}
Similarly, the fifteen specializations of $f_{36}(u,v,x)$ in the table all have Galois group 
$S_{36}$ except the bulleted one, for which the Galois group is intransitive.  
A presentation for this algebra $K_{36,-1/4,-1/4}$ is $\Q[x]/(g_{10}(x) g_{13}(x^2))$ with 
\begin{eqnarray*}
g_{10}(x) & = & x^{10}-4 x^9+2 x^8+8 x^7-8 x^6+8 x^5-20 x^4-10 x^2+80 x-60, \\
g_{13}(x) &=  &  x^{13}-44
    x^{12}+728 x^{11}-5256 x^{10}+15240 x^{9}-5320 x^{8}-41620 x^{7} \\ && +72280
    x^{6} -33940 x^{5}-4320 x^4-8760 x^3+20480 x^2-6140 x+480.
\end{eqnarray*}
The field $\Q[x]/g_{10}(x)$ has Galois group $S_{10}$ and discriminant  $2^{25}  3^{6} 5^5$ while
for $\Q[x]/g_{13}(x)$ these invariants are  $S_{13}$ and $2^{33}  3^{9} 5^{11}$.  Despite the small
exponents, these fields are wildly ramified not only at $2$, but also at $3$ and $5$.

\subsection{Larger degrees}
\label{largerdegrees}
In larger degrees, the numerics of the sets $U_\nu(\Z^P)$ are reflected more
clearly in the number fields constructed.  For example, in a degree $202$ example
of \cite{RHNF}, the specializations at $u \in U_{21^3}(\Z^{\{2,3,5\}})$ 
produce $2947$ distinct fields, all full in the sense of having Galois group all of $A_{202}$ 
or $S_{202}$, all wildly ramified at $2$, $3$, and $5$, and
unramified elsewhere.    We similarly expect $U_{21^3}(\Z^{\{2,3,5\}})$ to be
likewise responsible for exactly $2947$ distinct
full fields in many degrees $m>202$.   It
seems possible that the Hurwitz construction 
accounts for all full fields in $NF_m(\{2,3,5\})$
for most of these degrees $m$.    

As another example which gives a numerical sense of the asymptotics 
of this situation, 
consider the specialization set 
$U_{32768,1^3}(\Z^{\{2,3,5\}}) = \Polys_{[32768]}(\{2,3,5\})$,
chosen because it contains $s_{8,1}$ from 
Figure~\ref{fractalroots}.   From the
generating function \eqref{genfunct235}, this specialization set
contains more than $7.46 \times 10^{43}$ elements.  One
of the smallest degree covers 
of $U_{32678,1^3}$, in the language
of \cite{RHNF,RV}, comes from
the Hurwitz parameter $(S_5,(21^3,32,221,5),(32768,1,1,1))$.
This cover has degree exactly 
$(10^{32768} \cdot 20 \cdot 15 \cdot 24)/ (60 \cdot 120) = 10^{32768}$.
As $u$ ranges over the large set $U_{32768,1^3}(\Z^{\{2,3,5\}})$
the specialized algebras $K_u$ are all ramified
within $\{2,3,5\}$.    Other Hurwitz parameters give this
 same degree and we expect that there are
 many full fields in $NF_{10^{32768}}(\{2,3,5\})$.  
 The point for this paper is that polynomials
 with bad reduction within $\{2,3,5\}$ 
 are an ingredient in the construction of
 these $K_u$.   By way of contrast,
 it seems possible that $NF_{10^{32678}+1}(\{2,3,5\})$ is empty.

\section{Future directions}
\label{conclude}
       
\subsection{Specialization sets $U_\nu[\Z^P]$ for general $\nu$}  
 Let $\nu = (\nu_1,\dots,\nu_r)$ be a sequence of positive integers.  For 
$F$ a field, let $\Conf_\nu(F)$ be the set of tuples of disjoint divisors $(D_1,\dots,D_r)$ on
the projective line over $F$, with $D_i$ consisting of $\nu_i$ distinct geometric points.  The group $PGL_2(F)$ 
acts on $\Conf_\nu(F)$ by fractional linear transformations.  The object 
$\Conf_\nu$ itself is a scheme which is smooth and of
relative dimension $\sum \nu_i$ over $\mbox{Spec}(\Z)$.  

There is a natural quotient scheme $U_\nu = \Conf_\nu/PGL_2$.  The
map $\epsilon_\nu : \Conf_\nu \rightarrow U_\nu$ induces
a bijection $\Conf_\nu(F)/PGL_2(F) \rightarrow U_\nu(F)$ whenever
$F$ is an algebraically closed field.  
In the case that $\nu_{r-2} = \nu_{r-1}=\nu_{r}=1$,
the action of $PGL_2(F)$ on $\Conf_\nu(F)$ is free for 
all fields $F$, and the maps $\Conf_\nu(F)/PGL_2(F) \rightarrow U_\nu(F)$
 are always bijective.  
The general case is more complicated because
there may be points in $\Conf_\nu(F)$ for which the
stablizer in $PGL_2(F)$ is nontrivial.  
 The proofs of Theorems~\ref{thm2111} and 
\ref{thm3111} involved the maps $\epsilon_{\nu}$
for $\nu = 21^2$ and $\nu=31$ without 
using this notation.   Via the coordinates $w$ and $j$ respectively,
one has $U_{21^2}(F) = F^\times$ and $U_{31}(F) = F$ 
for $F$ of characteristic $>2$ and $>3$ respectively.  
The complications with fixed points are above $w=1$ 
and $j=0,1$.  

For $P$ a finite set of primes, let $U_\nu[\Z^P]$ be 
the image of $\Conf_\nu(\Z^P)$ in $U_\nu(\Q)$.  
The set $U_\nu[\Z^P]$ may be strictly smaller than
the set $U_\nu(\Z^P)$ of scheme-theoretical $\cP$-integral points, 
as illustrated by the equalities
$U_{21^2}[\Z^P] = T_{\infty,2,\infty}(\Z^P) \cup \{1\}$ 
and $U_{31}[\Z^P] = T_{3,2,\infty}(\Z^P) \cup \{0,1\}$,
which hold respectively under the 
assumptions $\{2\} \subseteq P$ and
$\{2,3\} \subseteq P$.  

In this paper, we have focused on tabulating 
$\Polys_\kappa(\cP)$ to keep sets small and
have a clear graph-theoretic interpretation.  
However from the point of view of 
Section~\ref{specialization}, our 
actual problem has been the identification
of $U_\nu(\Z^P)$ whenever $\nu_{r-2} = \nu_{r-1}=\nu_{r}=1$.
The natural generalization is to identify $U_\nu[\Z^P]$ 
for general $(\nu,P)$.   The Katz and Hurwitz 
theories of the previous section naturally 
give covers of general $U_\nu$ for general $\nu$.   

The general problem of identifying $U_\nu[\Z^P]$ has 
the same character as the special
case that we treat, but is technically
more complicated because 
elements of $\Conf_\nu(\Z^P)$ 
can no longer be canonically normalized
by applying a fractional linear transformation.
In the extreme case $\nu=(n)$
the complications become
quite severe:  describing the
scheme $U_n$ is a goal of
 classical invariant theory,
 and explicit results become 
 rapidly more complicated as $n$ increases.  
 
 The group $S_n$ acts naturally on the scheme $U_{1^n}$.
 Despite the normalization of three points
 to $0$, $1$, and $\infty$ in previous sections, 
 the influence of this automorphism group 
 has been visible.  For example, 
 the natural automorphism group of the left
 half of Figure~\ref{pictboth} is $S_5$, and it acts transitively
 on the twelve components of
  $U_{1^5}(\R)$. 
 An alternative viewpoint on $U_\nu$ for 
 general $\nu = (\nu_1,\dots,\nu_r)$ is
 via the equation
 \begin{equation}
 \label{Ucover}
 U_\nu = U_{1^n}/(S_{\nu_1} \times \cdots \times S_{\nu_r}).
\end{equation}
For example, the left half of Figure~\ref{pictboth} covers
the right half via $U_{1^5} \rightarrow U_{21^3}$, 
$(s,t) \mapsto (u,v) = (s+t,(1-s)(1-t))$.   
The map is not surjective even on $\R$-points 
or $\Q_p$-points.  The map is  
very far from far from surjective on the
 $\Z^P$-points of interest to us, and so the
 new $U_\nu$ present genuinely new 
 arithmetic sets $U_\nu[\Z^P]$ to
 be identified, despite
 the tight relation \eqref{Ucover}.        
 Birch and Merriman \cite{BM} proved, as part of a considerably larger theory, that 
the sets $U_\nu[\Z^P]$ are all finite.  Their finiteness theorem was made effective by
Evertse and Gy\H{o}ry \cite{EG}.   

\subsection{Descriptions of $U_\nu[\Z^P]$}  
The most straightforward continuation of this paper
would be to completely identify more $U_\nu[\Z^P]$.   
Staying first in the  context 
$\nu_{r-2}=\nu_{r-1}=\nu_r=1$, 
the direction needing most
attention is general completeness
results for the excellent $P$-units
introduced in Section~\ref{genpolys}.  
{\em Magma} \cite{Mag} already has the efficient command \verb@ExceptionalUnits@ giving 
complete lists of exceptional units.  
An extension of its functionality 
to $P$-units would immediately move many
$U_\nu(\Z^P)$ currently in the
second regime of expected completeness into
the first regime of proved completeness.  Leaving the context of 
$\nu_{r-2} = \nu_{r-1} = \nu_r=1$, there
are many more $(\nu,P)$ for which 
complete identification of $U_\nu[\Z^P]$ is
within reach, as it is only required that 
the normalization 
problem be resolved in some different
way.  

In the third regime of the introduction, where
complete tabulation is impossible, 
there are still many 
questions to pursue.  One would 
first like heuristic estimates 
on $|U_\nu[\Z^P]|$; the study
of exceptional units in \cite{Ni} 
looks to be a useful guide. 
The ``vertical'' direction of $\cP$ fixed and
$\nu$ varying is interesting
from the point of view of constructing
number fields with larger degree and
 bounded ramification.
In this direction
it seems that close attention 
to constructional 
techniques like those of 
Section~\ref{largedegree}
may yield good lower bounds.  
In the  ``horizontal'' direction of 
$\nu$ fixed and $\cP$ increasing,
the Reduction Bound~\ref{trivbound} becomes
particularly important and
upper bounds on $|U_\nu[\Z^P]|$
may be available.  
Finally, the $U_\nu[\Z^P]$ 
are not just bare sets to be 
tabulated or counted, 
one should also pay attention
to their natural structures.  
Figure~\ref{pictboth} suggests
that in the horizontal direction
the asymptotic distribution
of $U_\nu[\Z^P]$ 
in $U_\nu(\R)$ may
be governed by interesting
densities.    The asymptotic distribution
of $U_\nu[\Z^P]$ in $U_\nu(\Q_p)$ 
is important for understanding
the $p$-adic ramification of
number fields constructed
via covers, and may likewise
be governed by densities.


\begin{thebibliography}{10}
       
\bibitem{BR} Jos\'e Bertin and Matthieu Romagny.  
 Champs de Hurwitz. M\'em. Soc. Math. Fr. (N.S.) No. 125-126 (2011), 219 pp. 

\bibitem{BM} B.\ J.\ Birch and J.\ R.\ Merriman.
Finiteness theorems for binary forms with given discriminant. 
Proc. London Math. Soc. (3) 24 (1972), 385--394.

\bibitem{Mag} W. Bosma, J. J. Cannon, C. Fieker, A. Steel (eds.), Handbook of Magma Functions, Edition 2.19 (2012), 5478 pages.

\bibitem{Co} F.\ B.\ Coghlan.  Elliptic curves with condutor $N = 2^m 3^n$.  Ph.\ D.\ Thesis, Manchester University (1967).  Tables in: Modular Forms of One Variable IV. Springer Lecture Notes in Math. 476 (1975) pp. 123-134. Springer Verlag.



\bibitem{Cr} J. E. Cremona.  \verb+http://www.warwick.ac.uk/staff/J.E.Cremona/ftp/data/extra.html+



\bibitem{Jo} John W. Jones. Number fields unramified away from 2.  
J. Number Theory 130 (2010), no. 6, 1282--1291. 

\bibitem{Database} John W.\ Jones and David P. Roberts.  A database of number fields.  In preparation.  Database at \verb+http://hobbes.la.asu.edu/NFDB/+

\bibitem{EG}  J.\-H.\ Evertse and K.\ Gy\"ory. 
Effective finiteness results for binary forms with given discriminant. 
Compositio Math. 79 (1991), no. 2, 169Ð204. 

\bibitem{Ka} Nicolas Katz.  Rigid Local Systems.  Annals of Mathematics Studies, 139. Princeton University Press, Princeton, NJ, 1996. viii+223 pp.


\bibitem{LN}  Armin Leutbecher and Gerhard Niklash.  
On cliques of exceptional units and Lenstra's construction of Euclidean fields. Number theory (Ulm, 1987), 150--178, 
Lecture Notes in Math., 1380, Springer, New York, 1989. 

\bibitem{Ni} Gerhard Niklasch.
Counting exceptional units. 
Journ\'ees Arithm\'etiques (Barcelona, 1995). 
Collect. Math. 48 (1997), no. 1-2, 195--207. 

\bibitem{RABC} David P.\ Roberts.  An $ABC$ construction of number fields.
 Number theory, 237--267,
CRM Proc. Lecture Notes, 36, Amer. Math. Soc., Providence, RI, 2004.

\bibitem{RFract} \bysame Fractalized cyclotomic polynomials. 
Proc. Amer. Math. Soc. 135 (2007), no. 7, 1959--1967. 

\bibitem{RHNF} \bysame Hurwitz number fields.  In preparation.

\bibitem{RM05} \bysame Covers of $M_{0,5}$ and  number fields.  In preparation.

\bibitem{RV} David P.\ Roberts and Akshay Venkatesh.  Hurwitz monodromy and full number fields.   arXiv:1401.7379.  Submitted.

\bibitem{dW} B. M. M. de Weger.  Solving exponential Diophantine equations using lattice basis reduction algorithms.  
J. Number Theory 26 (1987), no. 3, 325--367.
\end{thebibliography}
\end{document}